\documentclass[openright,a4paper,11pt]{amsart}

\usepackage[latin1]{inputenc}
\usepackage[T1]{fontenc}
\usepackage{amsmath}
\usepackage{amssymb}

\usepackage[frenchb]{babel}

\usepackage{amsfonts}
\input xy
\xyoption{all}
\usepackage[dvips]{graphicx}
\usepackage{boxedminipage}
\usepackage{color}

\newcommand{\ZZ}{{\mathbb Z}}

\newcommand{\CC}{{\mathbb C}}
\newcommand{\RR}{{\mathbb R}}

\newcommand{\TT}{{\mathbb T}}
\newcommand{\A}{{\mathcal A}}
\newcommand{\td}{{\text{''}}}
\newcommand{\tg}{{\text{``}}}

\newtheorem{thm}{Théoreme}[section]
\newtheorem{defi}[thm]{Définition}
\newtheorem{prop}[thm]{Proposition}

          {\theoremstyle{remark}
}
          {\theoremstyle{definition}
}

\newenvironment{exo}{\it{}{}}

\begin{document}
\title{Un peu de géométrie tropicale}
\author{Erwan Brugallé}
\address{Université Pierre et Marie Curie,  Paris 6, 175 rue du Chevaleret, 75 013 Paris, France}
\email{brugalle@math.jussieu.fr}

\date{\today}

\maketitle

Quelles figures étranges aux propriétés mystérieuses  se
cachent derrière ce nom  énigmatique de \textit{géométrie 
  tropicale}? Sous les tropiques comme ailleurs, il est difficile de trouver
plus simple qu'une droite. Ce sera notre premier objet d'étude.

\begin{figure}[h]
\begin{center}
\begin{tabular}{ccccc}
\includegraphics[width=3cm, angle=0]{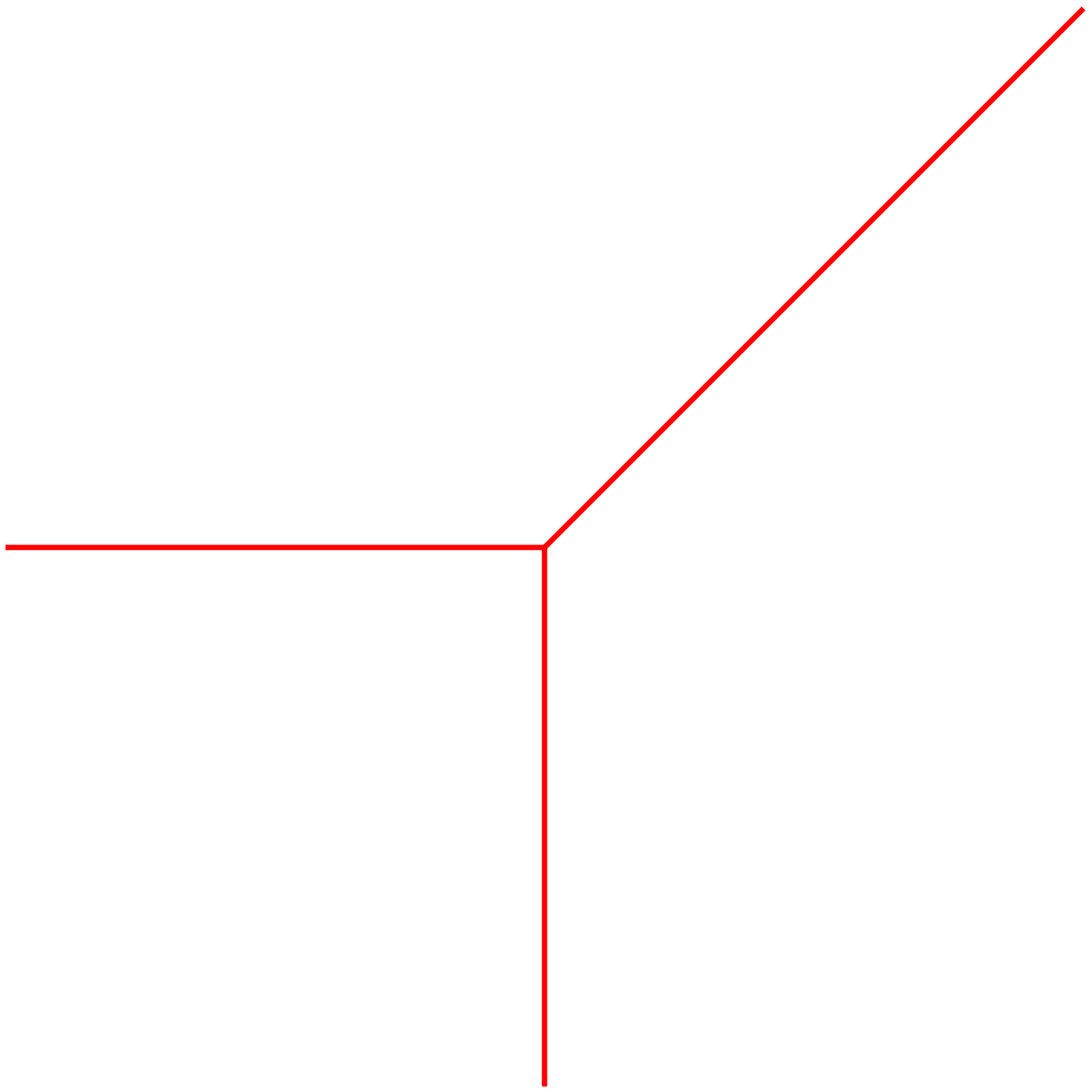}&\hspace{3ex} &
\includegraphics[width=3cm, angle=0]{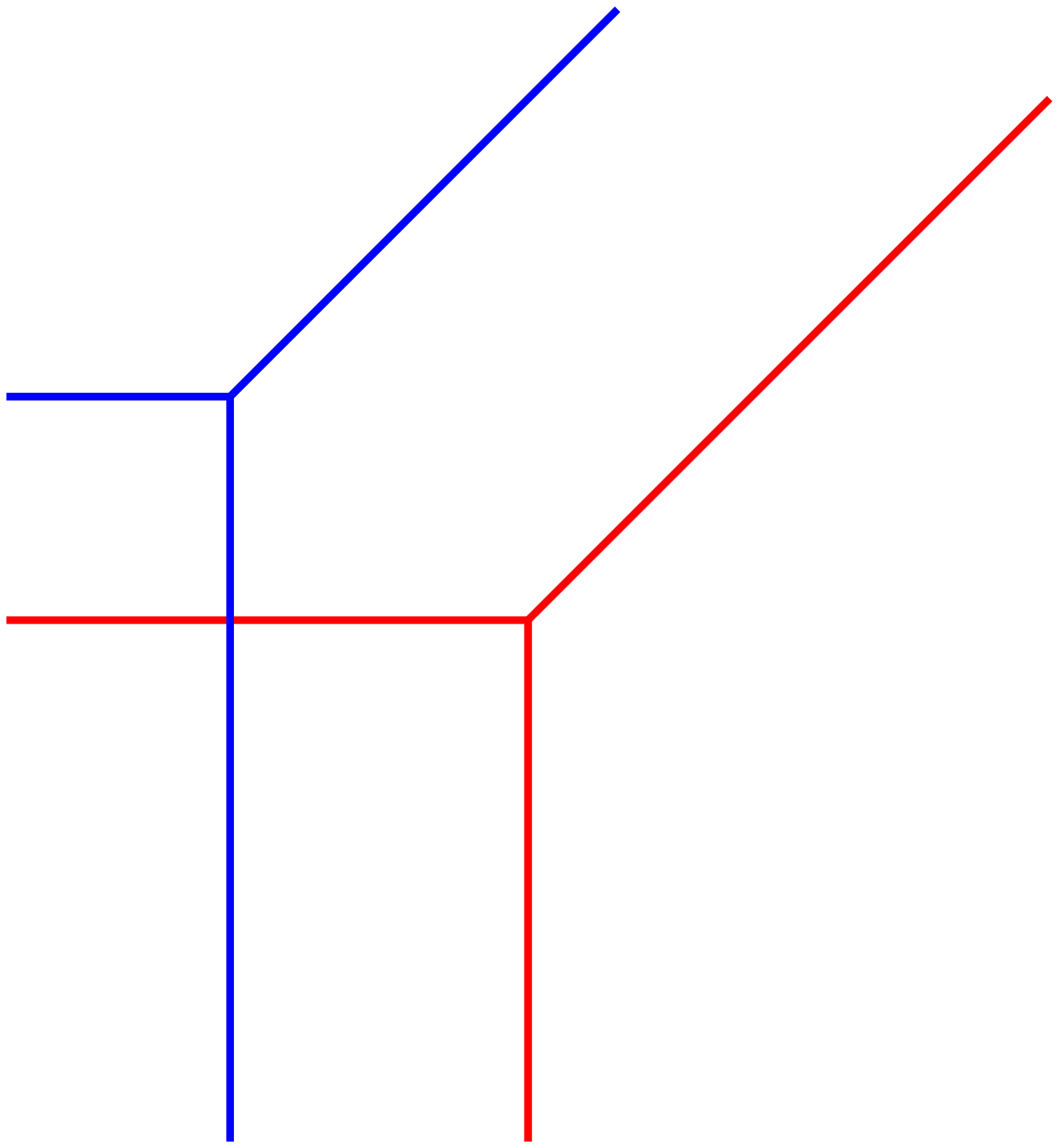}&\hspace{3ex} &
\includegraphics[width=3cm, angle=0]{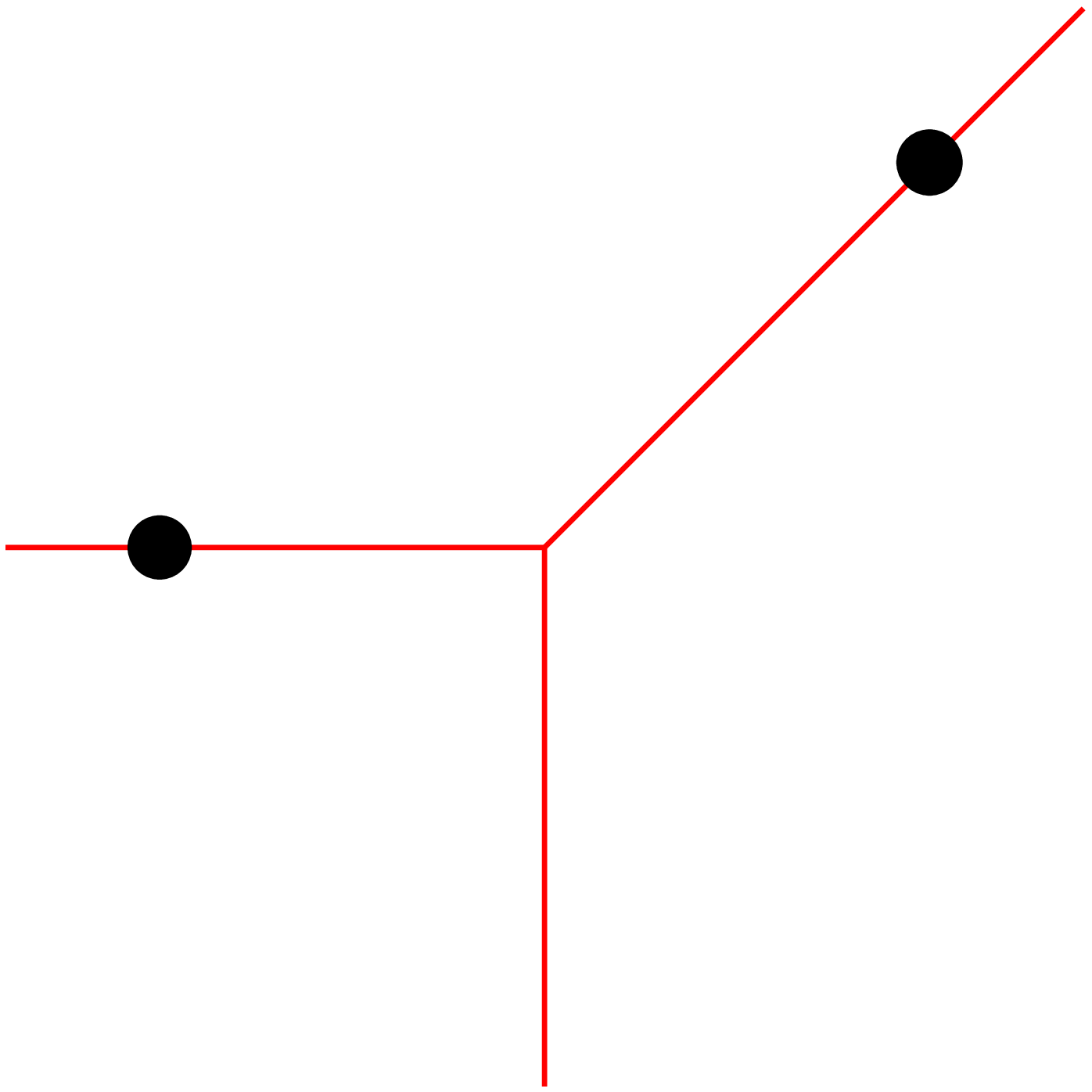}
\\ \\a)  && b)
 && c)
\end{tabular}
\end{center}
\caption{La droite tropicale}
\label{intro}
\end{figure}

  Une droite tropicale est
formée de 3 
demi-droites usuelles 
de directions $(-1,0)$, $(0,-1)$ et $(1,1)$
émanant d'un point quelconque du plan (voir figure \ref{intro}a). 
On peut se demander avec raison
pourquoi appeler une droite, tropicale ou
autre, 
cet objet bizarroïde...
En s'y penchant de plus près, on constate 
que ces droites tropicales satisfont les même propriétés
géométriques de base que les droites ``normales'' ou
 ``classiques'':  
deux droites tropicales  se coupent en un unique
point (voir figure \ref{intro}b), et deux points
 du plan définissent une unique droite tropicale (voir 
figure \ref{intro}c). 

Plus important encore, même si moins visible sur le
dessin,  droites classiques et tropicales sont toutes  deux
données par une équation de 
la forme $ax+by+c=0$. 
Dans le cadre de l'algèbre standard, c'est-à-dire en appelant une addition une addition et 
une multiplication une multiplication, on reconnaît sans peine une 
droite classique derrière cette équation. 
Mais dans le monde tropical,  additionner veut dire prendre le
  maximum,
 multiplier signifie additionner, et tous les objets changent de forme!
Pour dire toute la vérité, même ``être égal à 0'' prend un autre sens...

Géométries classique et tropicale sont donc élaborées suivant les mêmes
principes à partir de deux modes de calcul. Elles
 sont les visages de deux
algèbres différentes. 

\vspace{1ex}
La
géométrie tropicale n'est cependant pas qu'un jeu stérile pour mathématiciens
dés\oe uvrés. 
Le monde classique 
peut être  \textit{dégénéré} jusqu'au monde
tropical, et les objets tropicaux conservent alors naturellement certaines
propriétés des objets classiques dont ils sont  limite.
Ainsi, un énoncé tropical a de fortes chances d'avoir un énoncé
classique similaire. 
Or,  les objets tropicaux sont linéaires par morceaux et sont donc
beaucoup plus simples à étudier que leurs homologues classiques!

On pourrait donc résumer l'approche tropicale par le principe suivant: 
\begin{center}
\textit{Étudier des objets simples, énoncer des théorèmes sur des
objets compliqués}
\end{center}

\vspace{1ex}
Les premières parties de ce texte sont consacrées à l'algèbre tropicale,
 aux courbes tropicales et à quelques
 unes de
leurs propriétés. Nous expliquons ensuite pourquoi les géométries
classique et tropicale sont liées en montrant succinctement 
 comment
 le
monde classique dégénère jusqu'au  monde tropical. Puis
nous illustrons  le principe précédent  par  la
méthode dite du 
\textit{patchwork} pour construire des courbes algébriques réelles
 via les \textit{amibes}.
Nous terminons ce texte en donnant quelques références
bibliographiques.

\vspace{1ex}
Mais avant de rentrer dans le vif du sujet, il nous faut  expliquer
le pourquoi du mot  ``tropical''. Est-ce dû à la forme exotique des
objets considérés? À la présence d'amibes à squelette? 
 Avant de parler d'algèbre tropicale, on employait le
nom  plus prosaïque d'\textit{algèbre max-plus}. En l'honneur
des travaux de leur collègue brésilien Imre Simon, des chercheurs en
informatique de l'Université  Paris 7 décidèrent un jour de troquer
``max-plus'' pour ``tropical''. 
Laissons le mot de la fin à Wikipedia \footnote{15 mars 2009.} 
 sur
l'origine du mot 
``tropical'', \textit{it simply reflects the French 
  view on Brazil}.

\section{Algèbre tropicale}

\subsection{Opérations tropicales}
L'algèbre tropicale s'obtient en considérant l'ensemble des nombres
réels $\RR$ et en remplaçant l'addition par le maximum et la
multiplication par l'addition. En d'autres termes, on définit deux
nouvelles lois sur $\RR$, appelées \textit{addition} et \textit{multiplication
tropicales} et notées respectivement $\tg +\td$ et $\tg\times\td$, par
$$\tg x+y\td =\max(x,y) \ \ \ \ \ \ \ \ \tg x\times y\td =x+y $$
Dans tout ce texte, les opérations algébriques tropicales sont notées
entre guillemets. 
Comme pour la multiplication classique, nous abrégerons souvent $\tg
x\times y\td$ 
en $\tg x y\td$. 
Familiarisons nous avec ces deux opérations étranges en effectuant
quelques calculs  simples:
$$\tg 1+1\td=1, \ \ \tg 1+2\td=2, \ \ 
\tg 1+2+3\td=3, \ \  \tg 1\times 2\td=3,  \ \   \tg 1\times
(2+(-1))\td=3, $$
$$   \tg
1\times (-2)\td=-1, \ \ \tg (5+3)^2 \td= 10
$$

Ces deux lois tropicales ont beaucoup de  propriétés en commun avec
l'addition et la multiplication classiques. Par exemple, elles sont
toutes deux commutatives et la loi  $\tg
\times \td$ est distributive par rapport à la loi $\tg +\td$ (\textit{i.e.}
 $\tg (x+y)z\td =\tg xz+yz\td$). Il existe cependant deux
différences. Tout d'abord, l'addition tropicale n'a pas d'élément
neutre sur $\RR$. Qu'à cela ne tienne, nous pouvons
étendre naturellement nos deux 
opérations 
tropicales à $-\infty$ par 
$$\forall x\in\TT,  \ \ \ \ \tg x+(-\infty)\td =\max (x,-\infty) =x
 \ \ \ \ et \ \ \ \ \tg
x\times (-\infty)\td =x+(-\infty)=-\infty$$ 
où $\TT=\RR\cup\{-\infty\}$ est l'ensemble des \textit{nombres
  tropicaux}. Donc quitte à ajouter $-\infty$ à $\RR$, l'addition
tropicale a un élément neutre. 
En revanche il existe une
différence plus importante entre 
additions tropicale et 
classique: un élément de $\RR$ n'a pas 
de symétrique
pour la loi $\tg
+\td$. Autrement dit, il n'existe pas de soustraction tropicale.
De plus, cela ne marche pas cette fois   d'ajouter des
éléments à $\TT$ pour ``inventer'' des 
symétriques.
En effet $\tg +\td$ est  \textit{idempotente}, c'est-à-dire 
$\tg x+x\td=x $ pour tout $x$ dans $\TT$! Nous n'avons donc pas d'autre
choix que de nous accommoder de cette absence 
de symétriques
pour $\tg +\td$.

Mais mis à part ce dernier point,
l'ensemble $\TT$  muni des lois $\tg
+\td$ et $\tg \times\td$ satisfait toutes les autres propriétés d'un
corps. Par exemple $0$
est l'élément neutre de la multiplication tropicale, et tout élément
$x$ de $\TT$ différent de $-\infty$ a pour inverse $\tg\frac{1}{x}\td =-x$.
 On dit que $\TT$ est un \textit{semi-corps}. 

Attention  à ne pas aller trop vite dans l'écriture de
formules tropicales! Ainsi, $\tg 2x \td\ne \tg x+x\td$ mais $\tg 2x
\td=  x+2$, de même $\tg 1x\td\ne x$ mais  $\tg 1x\td= 
x+1$, ou encore $\tg
0x\td =x$ et $\tg (-1)x\td=x-1$.

\subsection{Polynômes tropicaux}
Après avoir défini l'addition et la multiplication tropicales,
nous arrivons naturellement à considérer des fonctions de la forme $P(x)=\tg
\sum_{i=0}^d a_ix^i\td$ avec  les $a_i$ dans $\TT$, c'est-à-dire des polynômes 
tropicaux\footnote{Nous considérons en fait 
  les fonctions 
  polynomiales plutôt que les polynômes.}. En réécrivant $P(x)$ avec
les notations classiques, on obtient 
$P(x)=\max_{i=1}^d(a_i + ix ) $. Voici quelques exemples de polynômes
tropicaux: 
$$\tg x\td = x, \ \ \  \tg 1+ x\td = \max(1,x), \ \ \  
\tg 1+ x +3x^2\td =\max(1,x,2x+3), $$
$$\tg 1+ x +3x^2+(-2)x^3\td =\max(1,x,2x+3, 3x-2) $$

Déterminons maintenant les racines d'un polynôme tropical. Mais
avant tout, qu'est ce qu'une racine tropicale?
Nous
rencontrons alors  un problème récurrent en mathématiques tropicales:
une notion classique a souvent plusieurs définitions équivalentes qui
ne le sont plus dans le monde tropical. Chaque définition  d'un même
objet classique produit
alors potentiellement autant d'objets tropicaux différents.

La première définition  d'une racine d'un polynôme  classique $P(x)$
est un élément $x_0$  tel que $P(x_0)=0$. Si on calque cette
définition en algèbre tropicale, on cherche alors
les éléments $x_0$ dans
$\TT$ tels que $P(x_0)=-\infty$. Or, si $a_0$ est le terme constant du
polynôme $P(x)$ alors  $P(x)\ge a_0$ pour tout $x$ dans $\TT$. Donc
si $a_0\ne - \infty$, le polynôme $P(x)$ n'a pas de racine... Cette
définition n'est pas très satisfaisante.

\begin{figure}[h]
\begin{center}
\begin{tabular}{ccccc}
\includegraphics[width=4cm, angle=0]{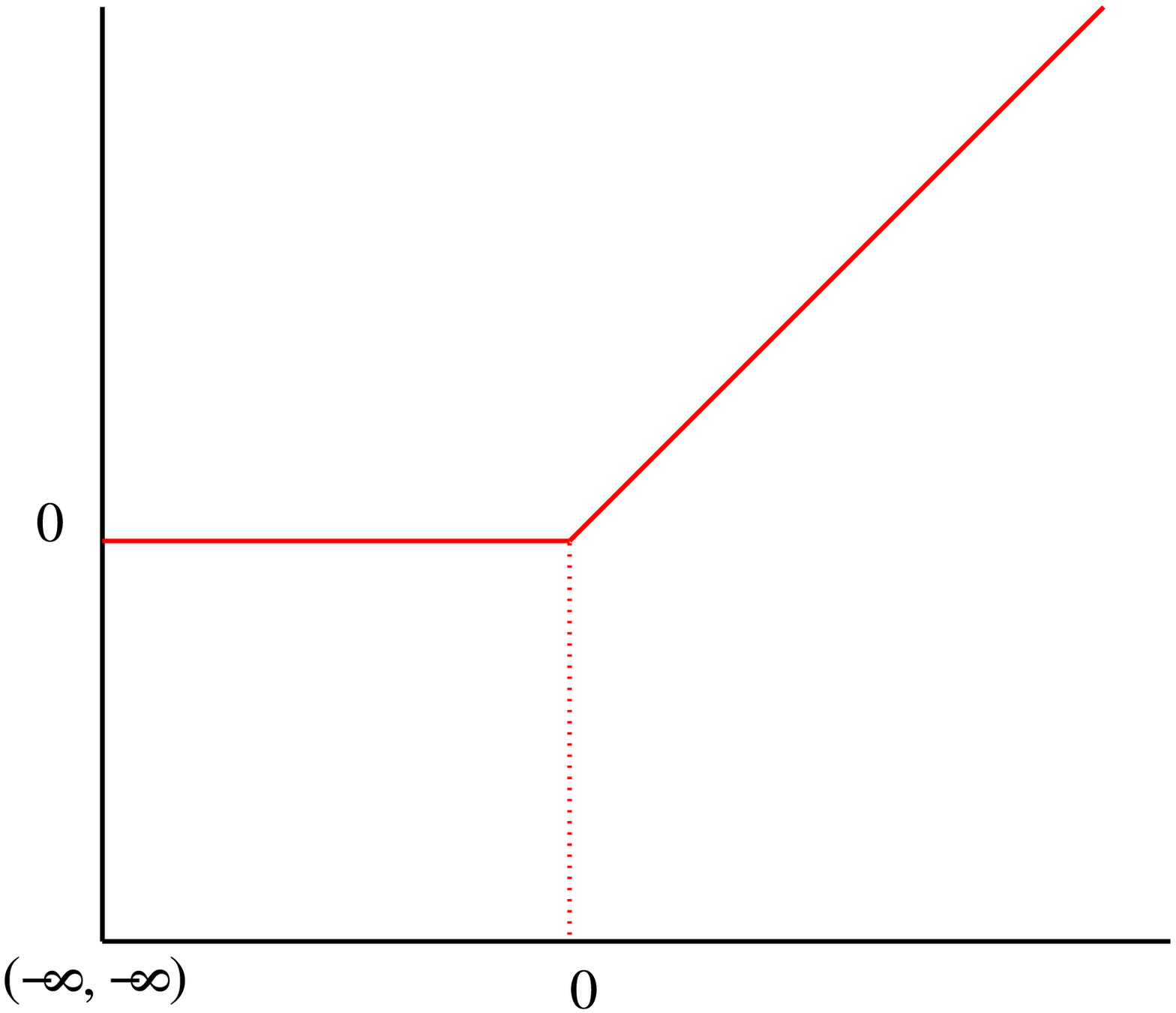}&\hspace{3ex} &
\includegraphics[width=4cm, angle=0]{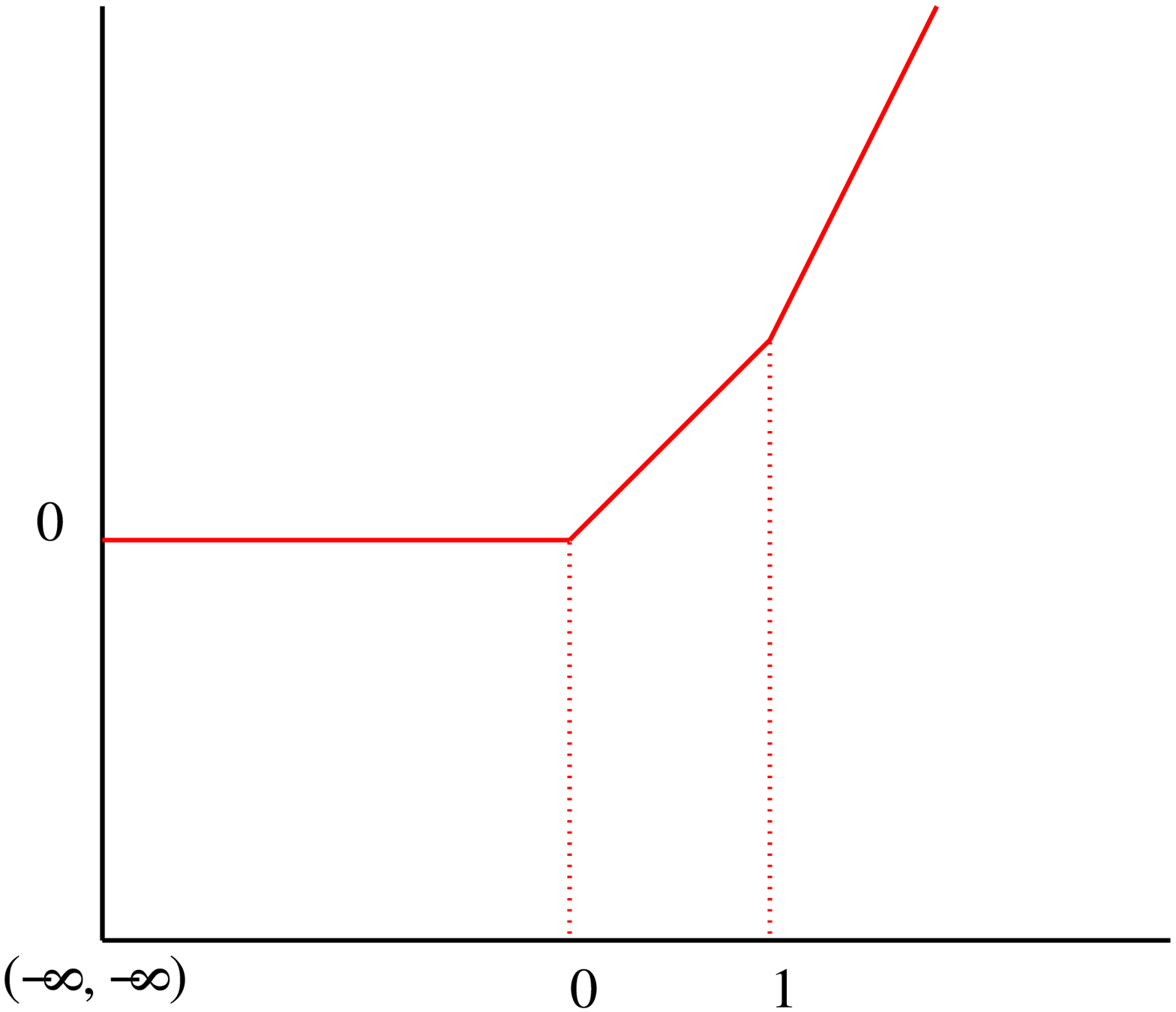}&\hspace{3ex} &
\includegraphics[width=4cm, angle=0]{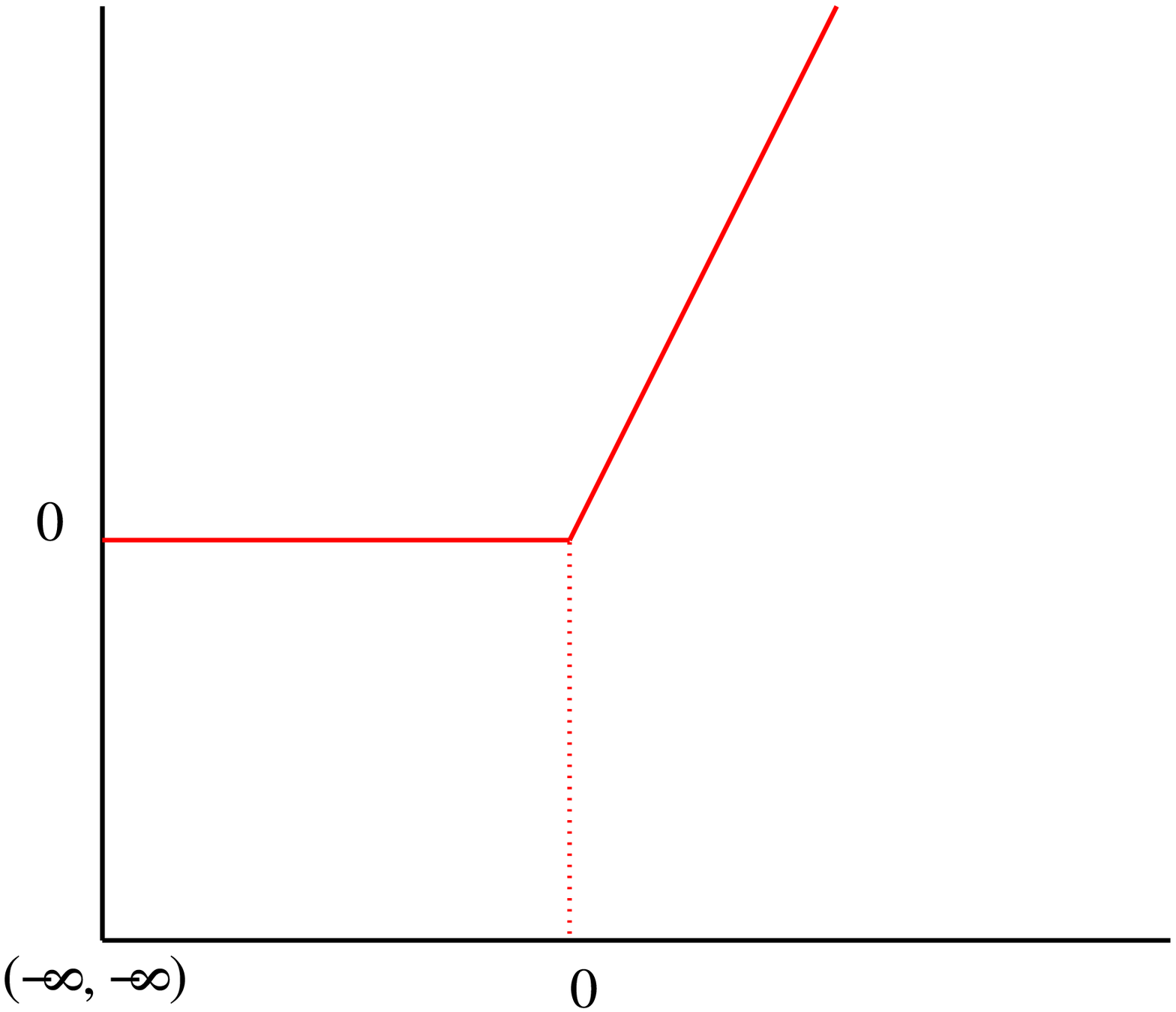}
\\ \\a) $P(x)=\tg 0+x \td$ && b) $P(x)=\tg 0+x + (-1)x^2\td$
 && b) $P(x)=\tg 0+x^2\td$
\end{tabular}
\end{center}
\caption{Exemples de graphe de polynômes tropicaux}
\label{graphes}
\end{figure}

Alternativement,  $x_0$ est une racine classique du polynôme
  $P(x)$ si il existe un polynôme $Q(x)$ tel que
$P(x)=(x-x_0)Q(x)$. Nous allons voir maintenant que cette définition est
la bonne en algèbre tropicale.
 Pour le comprendre, adoptons un point de vue géométrique sur le problème.
Un polynôme tropical est une fonction affine par morceaux (voir 
figure \ref{graphes}), et nous appelons \textit{racines tropicales} du
polynôme $P(x)$ tout point $x_0$ de $\TT$   pour lesquels le graphe de
$P(x)$ a un coin en $x_0$.
 De plus la
différence des deux pentes adjacentes à une racine tropicale
correspond à l'\textit{ordre} de cette racine.
 Ainsi, le polynôme $\tg 0+x\td$
a pour racine simple 0, le polynôme $\tg 0+x+ (-1)x^2\td$
a pour racines simples 0 et 1, et   le polynôme $\tg 0+x^2\td$
a pour racine double 0.

Les racines tropicales du polynôme $P(x)=\tg
\sum_{i=0}^d a_ix^i\td=\max_{i=1}^d(a_i + ix )$ sont donc
exactement les nombres tropicaux $x_0$ pour lesquels 
il existe $i\ne j$ tels que 
$P(x_0)=a_i+ix_0=a_j+jx_0$. On dit que le maximum de $P(x)$ est
atteint deux fois (au moins) en $x_0$.
Dans ce cas, l'ordre de $x_0$ est le maximum de
$|i-j|$ pour tous les 
$i$ et $j$ possibles qui réalisent ce maximum.
Par exemple, le maximum de $P(x)=\tg 0+x+x^2\td$ est atteint 3 fois en 0 et
l'ordre de cette racine est 2. 
De manière équivalente,  $x_0$ est une
racine tropicale d'ordre $k$ de $P(x)$ si 
il existe un polynôme $Q(x)$ tel que $P(x)=\tg (x+x_0)^kQ(x) \td$. 
Notons que le facteur $x-x_0$ en algèbre classique s'est transformé en le
facteur  $\tg x+x_0\td$, puisque
  la racine du polynôme
$\tg x+x_0\td$ est 
$x_0$ et non pas $-x_0$. 

Cette définition de racine tropicale a l'air nettement plus satisfaisante que
la première. De fait, on a la proposition suivante.

\begin{prop}
Le semi-corps tropical est algébriquement clos, c'est-à-dire que tout
polynôme tropical de degré $d$ a exactement $d$ racines tropicales comptées avec
leur multiplicité.
\end{prop}
Par exemple, on a les factorisations suivantes\footnote{Encore une
  fois, ces égalités 
  sont
vraies au niveau des fonctions polynomiales, pas au niveau des
polynômes!
Ainsi, $\tg 0+x^2\td$ et $ \tg (0+x)^2\td$ sont égales comme fonctions
polynomiales
mais pas comme polynômes.}:
$$\tg 0+x + (-1)x^2\td = \tg (-1)(x+0)(x + 1)\td \ \ \ et \ \ \ \tg 0 +
x^2\td = \tg (x+0)^2\td $$

\subsection{Exercices}
\begin{exo}
\begin{enumerate}
\item En quoi le fait que l'addition tropicale soit idempotente
  empêche-t-il l'existence de symétriques pour cette addition?
\item Tracer les graphes des polynômes tropicaux $P(x)=\tg x^3+2x^2+3x
  +(-1)\td$ et $Q(x)=\tg x^3+(-2)x^2+2x+(-1)\td$, et
  déterminer leurs racines tropicales.
\item Soit $a$ dans $\RR$ et $b$ et $c$ dans $\TT$. Déterminer les racines des
  polynômes tropicaux $ \tg ax+b\td$ et $ \tg ax^2+bx+c\td$. 
\end{enumerate}
\end{exo}

\section{Courbes tropicales}

\subsection{Définition}
N'ayons peur de rien, augmentons  le nombre de variables
de nos polynômes. Un polynôme tropical en deux variables s'écrit
$P(x,y)=\tg \sum_{i,j}a_{i,j}x^iy^j\td $, ou encore
$P(x,y)=\max_{i,j}(a_{i,j}+ix+jy)$ en notations classiques.
Ainsi $P(x,y)$ est encore une
fonction affine par morceaux, et  la \textit{courbe tropicale} $C$ définie
par $P(x,y)$ est  le lieu des coins de cette fonction.
 Autrement dit, $C$ est constituée des
 points $(x_0,y_0)$ de $\TT^2$ 
 pour lesquels le maximum de  $P(x,y)$ est atteint au moins deux fois
 en  $(x_0,y_0)$.

Avouons tout de suite que 
 nous nous contenterons dans ce texte
d'étudier les 
courbes tropicales dans $\RR^2$ au lieu de $\TT^2$. Cela n'entame en
rien la généralité de ce qui est discuté  ici, en revanche
les définitions, les énoncés et les dessins deviennent plus simples et
compréhensibles. 

Regardons la droite tropicale
définie par le polynôme 
$P(x,y)=\tg \frac{1}{2}+2x +(-5)y \td$. Il faut donc chercher les  points
$(x_0,y_0)$ dans $\RR^2$ qui 
vérifient l'un des trois systèmes suivant:
$$2+x_0=\frac{1}{2}\ge -5+y_0, \ \ \ \ \ \ \ \ \ -5+y_0=\frac{1}{2}\ge 2+x_0,
\ \ \ \ \ \ \ \ \ 2+x_0=-5+y_0\ge \frac{1}{2} $$ 

\begin{figure}[h]
\begin{center}
\begin{tabular}{ccc}
\includegraphics[width=4cm, angle=0]{Figures/Droite.eps}&
\includegraphics[width=4cm, angle=0]{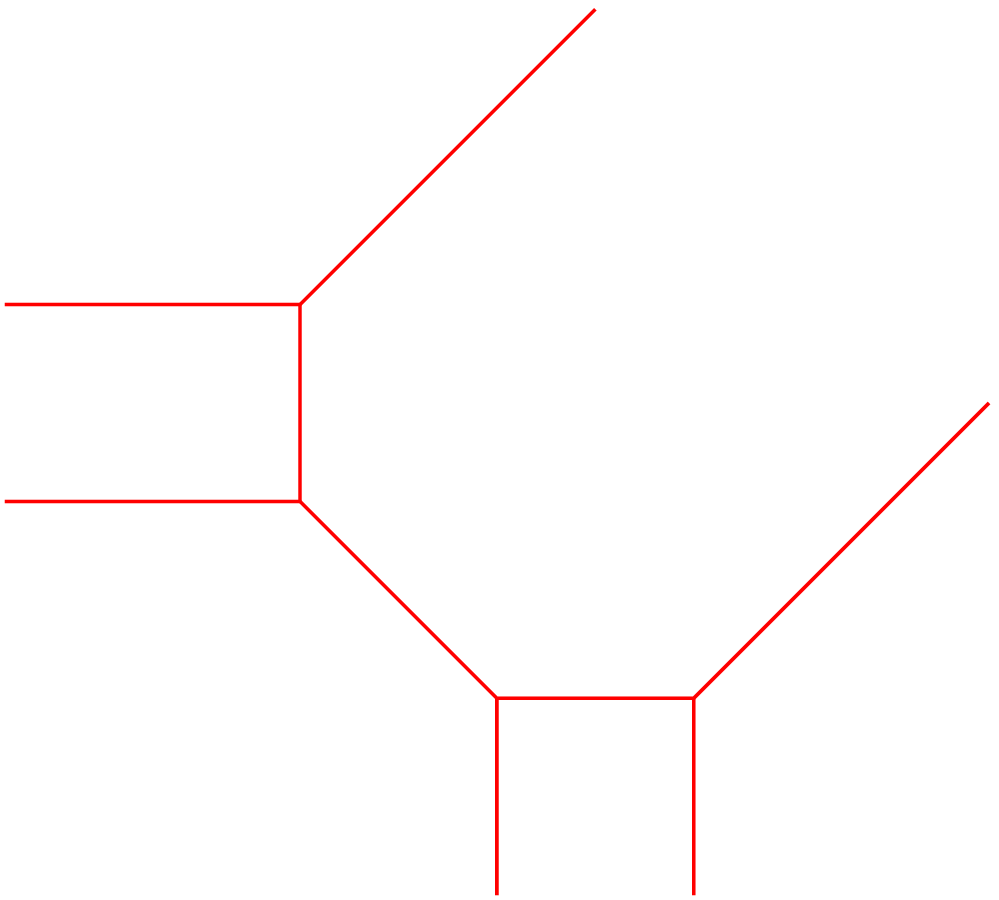}&
\includegraphics[width=4cm, angle=0]{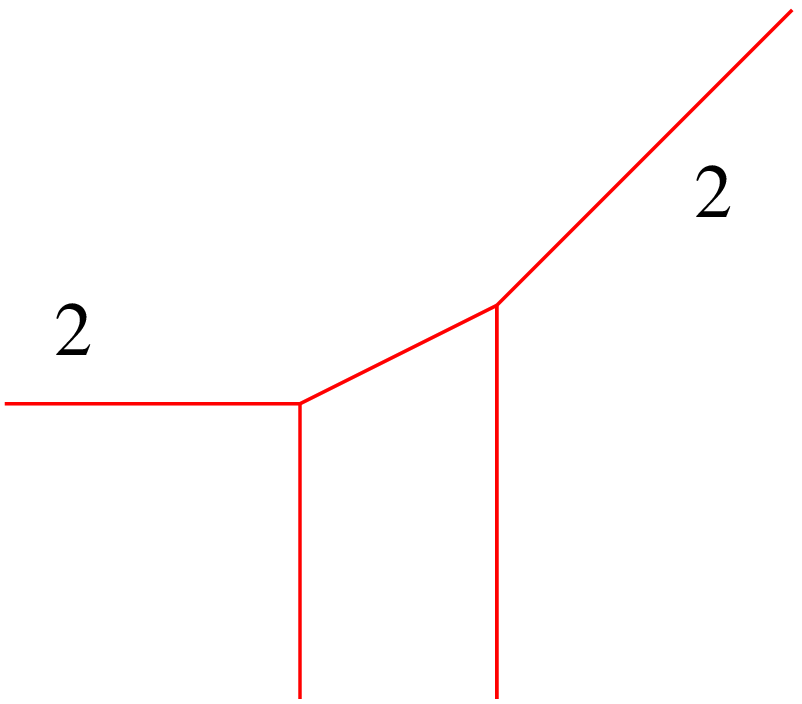}
\\ \\a) $\tg \frac{1}{2}+2x +(-5)y\td$ & b)  $\tg 3+ 2x + 2y + 3xy+y^2+x^2\td$ 
 & c) $\tg 0+ x +y^2+(-1)x^2\td$

\end{tabular}
\end{center}
\caption{Quelques courbes tropicales}
\label{droite}
\end{figure}

Notre droite tropicale est donc constituée des trois demi-droites
$\{(-\frac{3}{2},y) \ | \ y\le \frac{11}{2} \}$, $\{(x,\frac{11}{2}) \ |
\ x\le -\frac{3}{2} 
\}$, et  $\{(x,x+7) \ | 
\ x\ge  -\frac{3}{2} \}$ (voir  la figure \ref{droite}a). 

Il nous manque encore une donnée pour définir rigoureusement une courbe
tropicale. Le lieu des coins d'un polynôme tropical de deux variables est
constitué de segments et de demi-droites,  appelés
\textit{arêtes}, qui s'intersectent en des points, appelés
\textit{sommets}. 
Comme dans le cas des polynômes en une variable,    nous devons
prendre en compte  pour chaque
arête  la différence de pente de $P(x,y)$ des
deux côtés de cette arête. 
On arrive ainsi à la définition formelle
suivante.

\begin{defi}
Soit $P(x,y)=\tg \sum_{i,j} a_{i,j}x^iy^j\td $ un polynôme tropical. La courbe
tropicale $C$ définie  par $P(x,y)$  est l'ensemble des points $(x_0,y_0)$ de
$\RR^2$ tels qu'il existe $(i,j)\ne (k,l)$ vérifiant
$P(x_0,y_0)=a_{i,j}+ix_0+jy_0=a_{k,l}+kx_0+ly_0$.

Étant donné une arête de $C$, on définit le poids
de cette arête comme le maximum des $pgcd$ des nombres $|i-k|$ et
$|j-l|$ pour tous les couples $(i,j)$ et $(k,l)$ correspondant à cette arête.
\end{defi}
Dans les dessins, on n'écrira le poids d'une arête près de celle-ci
seulement si le poids est au moins 2. 
Dans le cas de la droite tropicale, toutes les arêtes sont de poids 1,
donc la figure \ref{droite}a représente bien une droite
tropicale. Deux exemples de courbes tropicales de degré 2 sont
représentés sur les figures \ref{droite}b et c. La conique tropicale de la
figure \ref{droite}c a deux arêtes de poids 2.

\subsection{Subdivision duale}
Un polynôme tropical $P(x,y)$ est donc donné par le maximum d'un nombre
fini de fonctions affines qui sont les monômes de $P(x,y)$.  De plus les
points du plan pour lesquels au moins deux des 
 monômes réalisent ce maximum sont précisément les
points de la courbe tropicale $C$ définie par $P(x,y)$. Affinons un
peu cette étude et pour chaque point $(x_0,y_0)$ de $C$,
considérons \textit{tous} les monômes de $P(x,y)$ qui réalisent ce
maximum en $(x_0,y_0)$.

\begin{figure}[h]
\begin{center}
\begin{tabular}{ccccc}
\includegraphics[width=1cm, angle=0]{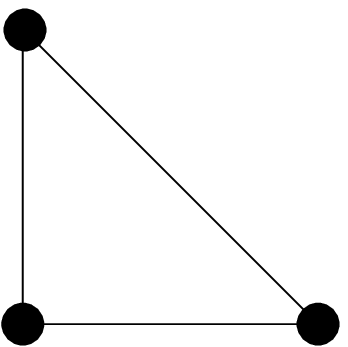}&\hspace{6ex} &
\includegraphics[width=2cm, angle=0]{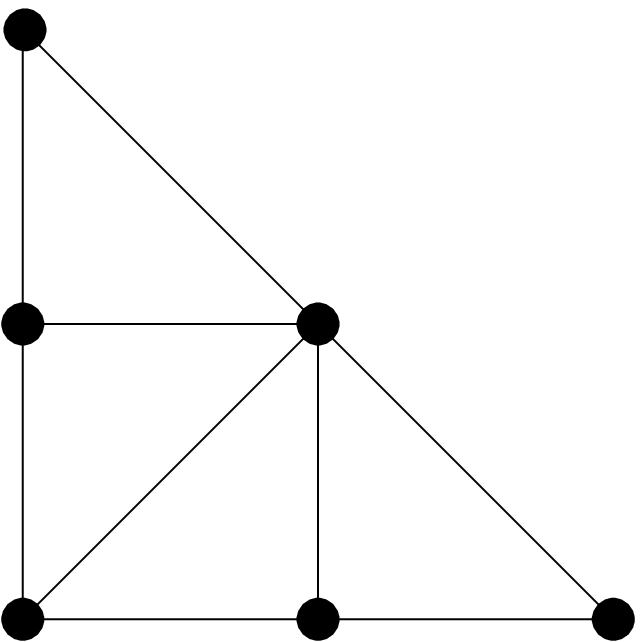}&\hspace{6ex} &
\includegraphics[width=2cm, angle=0]{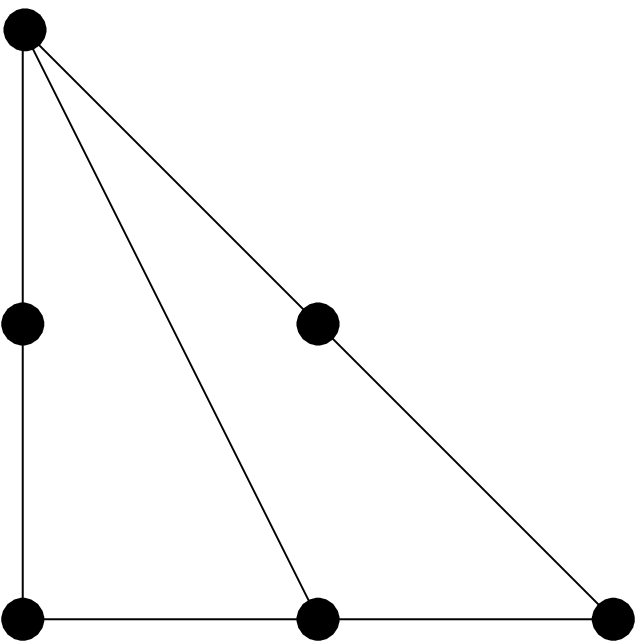} 
\\ \\a) && b) &&c)

\end{tabular}
\end{center}
\caption{Quelques subdivisions duales}
\label{subd}
\end{figure}

Étudions tout  d'abord  le cas de la droite tropicale $C$ définie par
$P(x,y)=\tg 
\frac{1}{2}+2x+(-5)y\td$ (voir 
figure \ref{droite}a). Au 
point $(-\frac{3}{2},\frac{11}{2})$, le sommet de la droite, les trois
monômes $\frac{1}{2}=\frac{1}{2}x^0y^0$, 
$2x=2x^1y^0$ et 
$(-5)y=(-5)x^0y^1$ ont la même valeur. 
Les exposants de ces monômes, c'est-à-dire
les points $(0,0)$,
$(1,0)$ et $(0,1)$, définissent un triangle  $\Delta_1$
(voir figure \ref{subd}a). Le long  de
l'arête horizontale de $C$, la valeur du polynôme $P(x,y)$ est donnée par les
monômes $0$ et $y$, c'est-à-dire par les monômes d'exposants $(0,0)$
et $(0,1)$. Le segment défini par ces deux exposants est donc l'arête
verticale du triangle $\Delta_1$. De même, les monômes donnant la
valeur de $P(x,y)$ le long  de l'arête
verticale (respectivement de pente $1$) de $C$ sont d'exposants $(0,0)$ et
$(1,0)$ (respectivement $(1,0)$ et $(0,1)$) et le segment défini par
 ces exposants est  l'arête horizontale (respectivement de pente
$-1$) du triangle $\Delta_1$.  

Que retenir de ce petit exercice? En
regardant  les monômes  donnant la valeur du polynôme tropical  $P(x,y)$ en un
point de la droite tropicale $C$, on s'aperçoit que le sommet de $C$
correspond au triangle $\Delta_1$ et que chaque arête $e$ de $C$
correspond à 
une arête de $\Delta_1$ dont la direction est perpendiculaire à celle
de $e$.
\vspace{1ex}

Regardons maintenant la
conique tropicale définie par le polynôme $P(x,y)=\tg
3+2x+2y+3xy+x^2+y^2 \td$  représentée à la figure \ref{droite}b. Cette
courbe a pour sommets les quatre 
points $(-1,1)$, $(-1,2)$, $(1,-1)$ et $(2,-1)$. En chacun de ces
sommets $(x_0,y_0)$, la valeur du polynôme $P(x,y)$ est donnée par trois monômes:
$$P(-1,1)=3=y_0+2=x_0+y_0+3
\ \ \ \ \ \ \ \ \ \ \ \ \ \ \ \ P(-1,2)=y_0+2=x_0+y_0+3=2y_0 \ \  $$
$$P(1,-1)=3=x_0+2=x_0+y_0+3
\ \ \ \ \ \ \ \ \ \ \ \ \ \ \ \ P(2,-1)=x_0+2=x_0+y_0+3=2x_0 \ \  $$
Ainsi pour chaque sommet de $C$, les
exposants des trois monômes correspondant définissent un triangle, et
ces quatre triangles sont
disposés comme sur la figure \ref{subd}b. De plus comme dans le cas de
la droite, 
pour chaque arête $e$ de $C$, les exposants des monômes
donnant la valeur de $P(x,y)$ le long de $e$ définissent une arête
d'un (ou deux) de ces triangles, et 
 la direction de cette arête est perpendiculaire à celle de $e$.

\vspace{1ex}
Expliquons maintenant ce phénomène en toute généralité. Soit
$P(x,y)=\tg \sum_{i,j} a_{i,j}x^iy^j\td$ 
un polynôme tropical quelconque. Le \textit{degré} de $P(x,y)$ est   le
maximum des sommes $i+j$ pour les coefficients $a_{i,j}$ différents de $-\infty$. Par simplicité, 
tous les polynômes de degré $d$ considérés dans ce texte vérifient
 $a_{0,0}\ne -\infty$, $a_{d,0}\ne  -\infty$ et  $a_{0,d}\ne
-\infty$. Ainsi, l'enveloppe convexe des points $(i,j)$ tels que
$a_{i,j}\ne -\infty$ est le triangle   $\Delta_d$ de
 sommets $(0,0)$, 
 $(d,0)$ et $(0,d)$.

Si $v=(x_0,y_0)$ est un sommet de $C$, alors l'enveloppe convexe des points
$(i,j)$ dans $\Delta_d\cap\ZZ^2$ tels que $P(x_0,y_0)=a_{i,j}+ix_0+jy_0$ est un polygone
$\Delta_v$ 
inclus dans $\Delta_d$. De même, si $(x_0,y_0)$ est un point à l'intérieur
d'une arête $e$ de $C$, alors l'enveloppe convexe des points
$(i,j)$ dans $\Delta_d\cap\ZZ^2$ tels que
$P(x_0,y_0)=a_{i,j}+ix_0+jy_0$ est un segment 
$\delta_e$ inclus dans $\Delta_d$.
Comme le polynôme tropical $P(x,y)$ est une fonction convexe affine
par morceaux, l'union des $\Delta_v$ forme une \textit{subdivision} de
$\Delta_d$. En d'autres termes l'union des polygones $\Delta_v$ est
égale au triangle $\Delta_d$, et  deux polygones $\Delta_v$ et
$\Delta_{v'}$ ont soit une arête commune, soit un sommet commun, soit ne s'intersectent pas.
 De plus, si $e$ est une arête de $C$ adjacente au sommet $v$,
alors $\delta_e$ est une arête de $\Delta_v$, et les directions
de $e$ et de $\delta_e$ sont perpendiculaires. Cette subdivision de
$\Delta_d$ est appelée la \textit{subdivision duale} à $C$.

Par exemple, les subdivisions duales aux courbes tropicales de la
figure \ref{droite} sont dessinées à
 la figure \ref{subd} (les points noirs
 représentent 
 les points à coordonnées entières, et ne sont pas nécessairement des
 sommets de la subdivision).

Remarquons que  $e$ est une arête de poids $w$ de $C$ si et
seulement si le
segment $\delta_e$  contient
$w+1$ points dans $\ZZ^2$. Ainsi, le degré d'une courbe tropicale peut
se lire directement sur la courbe: c'est la somme des poids des
arêtes infinies dans la direction $(-1,0)$ (ou $(0,-1)$, ou encore
$(1,1)$).
De plus, une courbe tropicale est donnée par sa subdivision duale,
à translation et longueur des arêtes
près.

\subsection{Graphes équilibrés et courbes tropicales}
\begin{figure}[h]
\begin{center}
\begin{tabular}{ccc}
\includegraphics[width=4cm,
  angle=0]{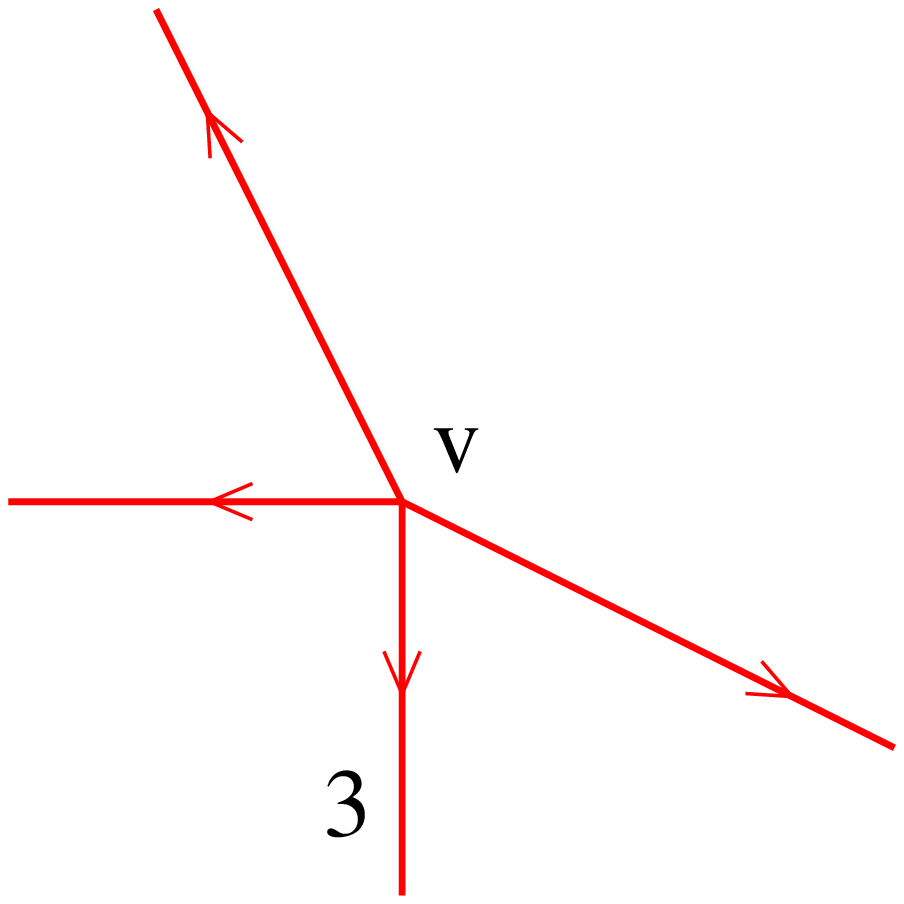}& \hspace{4ex} &
\includegraphics[width=4cm, angle=0]{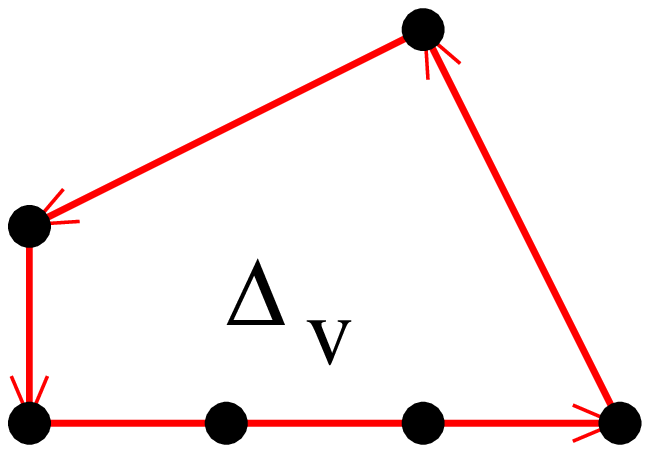}
\\ a) && b)
\end{tabular}
\end{center}
\caption{Relation d'équilibre}
\label{equ}
\end{figure}

La première conséquence de cette dualité est qu'une certaine équation,
appelée \textit{relation 
  d'équilibre}, est vérifiée en chacun des sommets d'une
courbe tropicale.
Soit $v$  un sommet de $C$ adjacent aux
arêtes $e_1,\ldots,e_k$ de poids respectivement
$w_1,\ldots,w_k$. Comme $e_i$ est 
supportée par une droite (au sens usuel) d'équation à coefficients
entiers, il existe un 
unique vecteur entier $\vec v_{i}=(\alpha,\beta)$ sur $e_i$ avec
$pgcd(\alpha,\beta)=1$ et d'origine le sommet $v$ (voir
figure \ref{equ}a). D'après la partie précédente, le polygone
$\Delta_v$ dual à $v$ se déduit immédiatement des vecteurs $w_1\vec
v_{1},\ldots, w_k\vec v_{k}$: si on oriente le bord de $\Delta_v$
dans le sens inverse des aiguilles d'une montre, alors  chaque arête
$\delta_{e_i}$ de 
$\Delta_v$ duale à
$e_i$ est obtenue à partir  du vecteur
$w_i\vec v_{i}$ par une rotation d'angle $\frac{\pi}{2}$ (voir figure
\ref{equ}b). 

Puisque $\Delta_v$ est fermé, 
 on lit alors immédiatement sur 
ce dernier la relation d'équilibre  suivante:
$$\sum_{i=1}^k w_i\vec v_i=0$$

Un graphe dans $\RR^2$ qui vérifie la relation  d'équilibre en chacun
des ces sommets est appelé un  \textit{graphe équilibré}. Nous venons
donc de voir que toute courbe tropicale est un graphe équilibré.
Il se trouve que la réciproque est vraie.

\begin{thm} 
Les courbes tropicales dans $\RR^2$ sont exactement les graphes
équilibrés.
\end{thm}

\begin{figure}[h]
\begin{center}
\begin{tabular}{ccccc}
\includegraphics[width=4cm, angle=0]{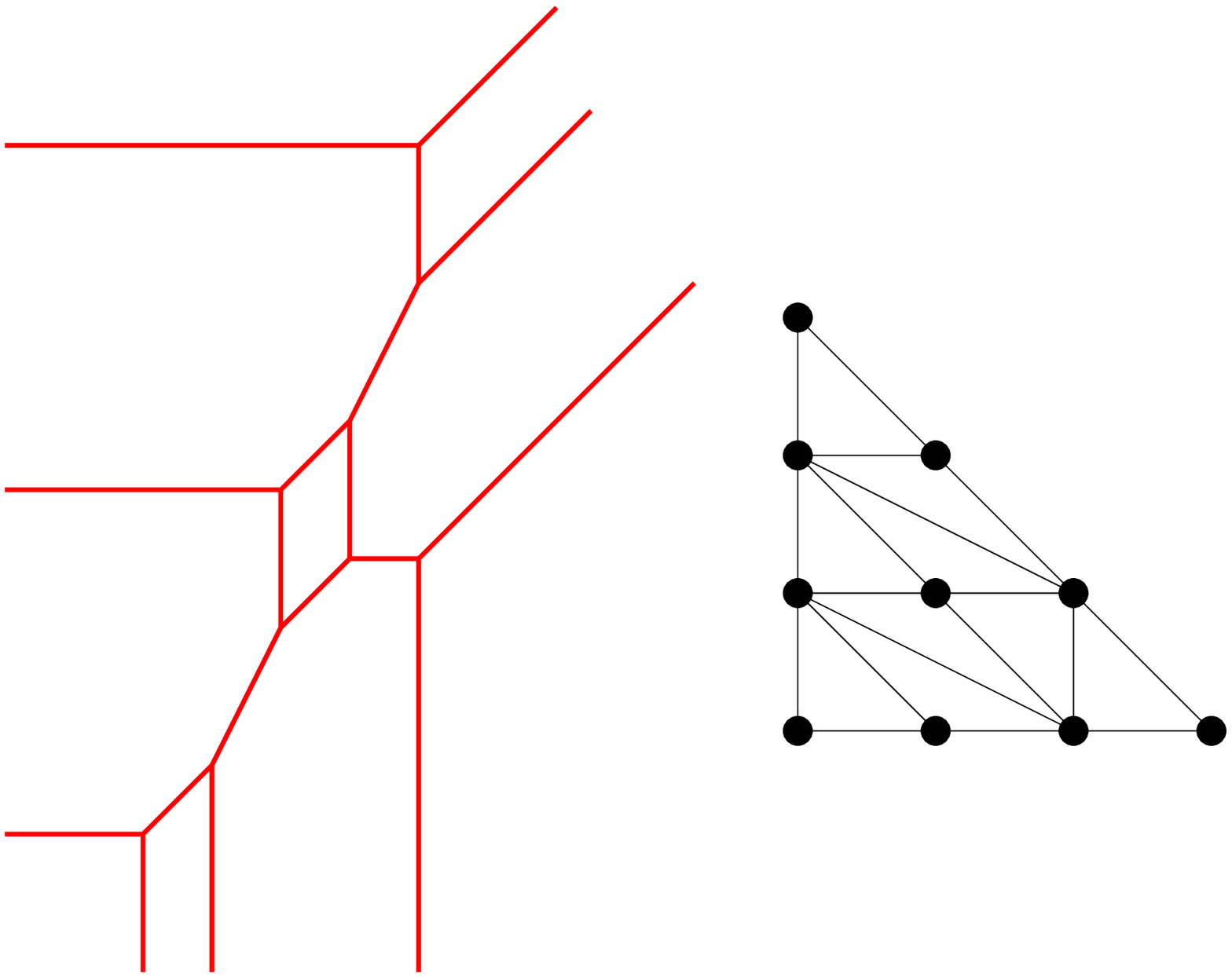}&\hspace{3ex} &
\includegraphics[width=4cm, angle=0]{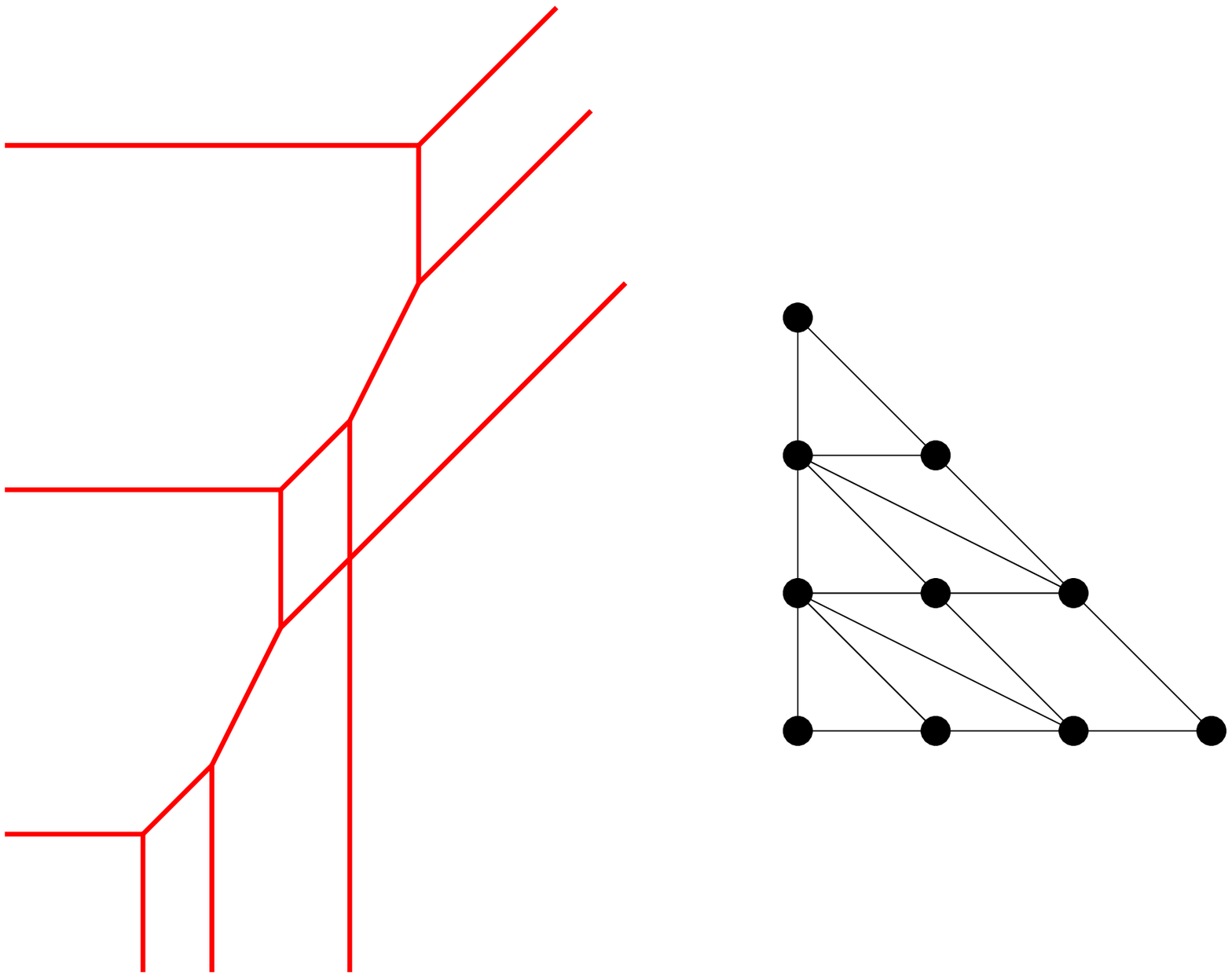}&\hspace{3ex} &
\includegraphics[width=4cm, angle=0]{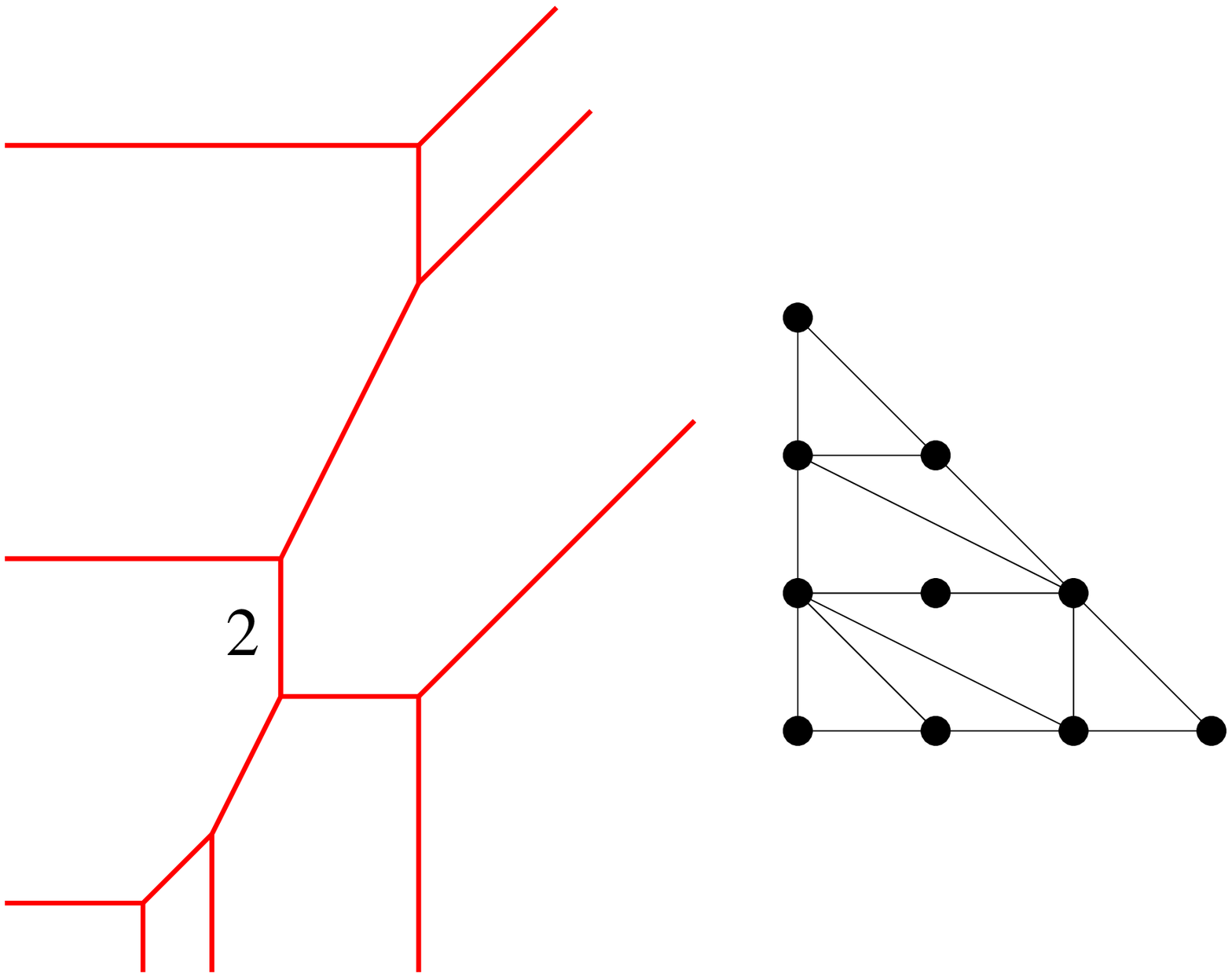}
\\ \\a) && b) &&c)

\end{tabular}
\end{center}
\caption{}
\label{equil}
\end{figure}

Ainsi, on peut affirmer qu'il existe des
polynômes tropicaux de degré 3 dont les courbes tropicales sont les graphes
équilibrés représentés à la figure \ref{equil}. Nous avons aussi
représenté dans chaque cas la subdivision  de $\Delta_3$ duale à la courbe.

\subsection{Exercices}
\begin{exo}
\begin{enumerate}
\item Dessiner les courbes tropicales définies par les polynômes tropicaux
  $P(x,y)=\tg 5 + 5x + 5y + 4xy+1y^2+x^2\td$ et $Q(x,y)=\tg  7+ 4x + y
  + 4xy+3y^2+(-3)x^2\td$, ainsi que  leur 
  subdivision duale.
\item  Un triangle tropical est un domaine borné de $\RR^2$ délimité
  par trois droites 
  tropicales. 
   Quelles sont les formes
  possibles d'un triangle tropical?
\item Montrer qu'une courbe tropicale de degré $d$ a au plus $d^2$ sommets.
\item Trouver une équation pour chacune des courbes  tropicales de la
  figure \ref{equil}. Le  rappel
  suivant peut être
  utile: si $v$ est un sommet d'une courbe tropicale définie par un 
polynôme tropical $P(x,y)$, alors  la valeur de $P(x,y)$ au voisinage
de $v$ est
donnée   uniquement par les monômes correspondant au
polygone dual à $v$.
\end{enumerate}
\end{exo}

\section{Intersection tropicale}

\subsection{Théorème de Bézout}
Un grand intérêt de la géométrie tropicale est de fournir un
modèle simple de la géométrie algébrique. Par exemple, 
les théorèmes
de base sur l'intersection de courbes tropicales nécessitent un bagage
algébrique nettement moins important que leurs homologues classiques. 
Nous allons
illustrer ce principe avec le théorème de Bézout qui affirme que deux
courbes algébriques  planes de degré $d_1$ et $d_2$ se coupent en
$d_1d_2$ points\footnote{Attention, ceci est
  un théorème 
  de géométrie projective! Par exemple, deux droites affines peuvent être
parallèles...}. Avant le cas général, regardons d'abord les droites et 
les coniques tropicales.

\begin{figure}[h]
\begin{center}
\begin{tabular}{ccccc}
\includegraphics[width=3cm, angle=0]{Figures/InterDte.eps}&\hspace{3ex} &
\includegraphics[width=4cm, angle=0]{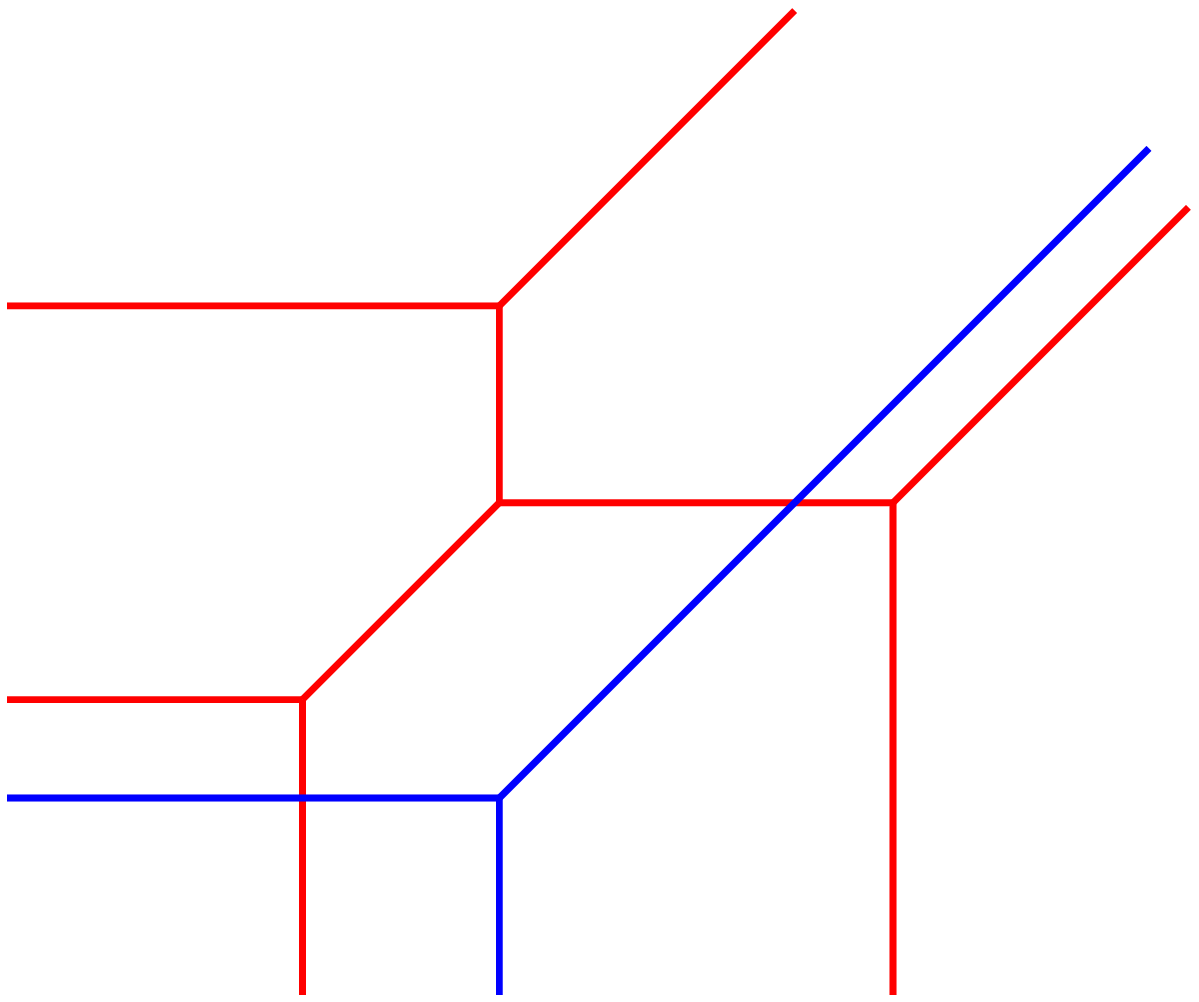}&\hspace{3ex} &
\includegraphics[width=4cm, angle=0]{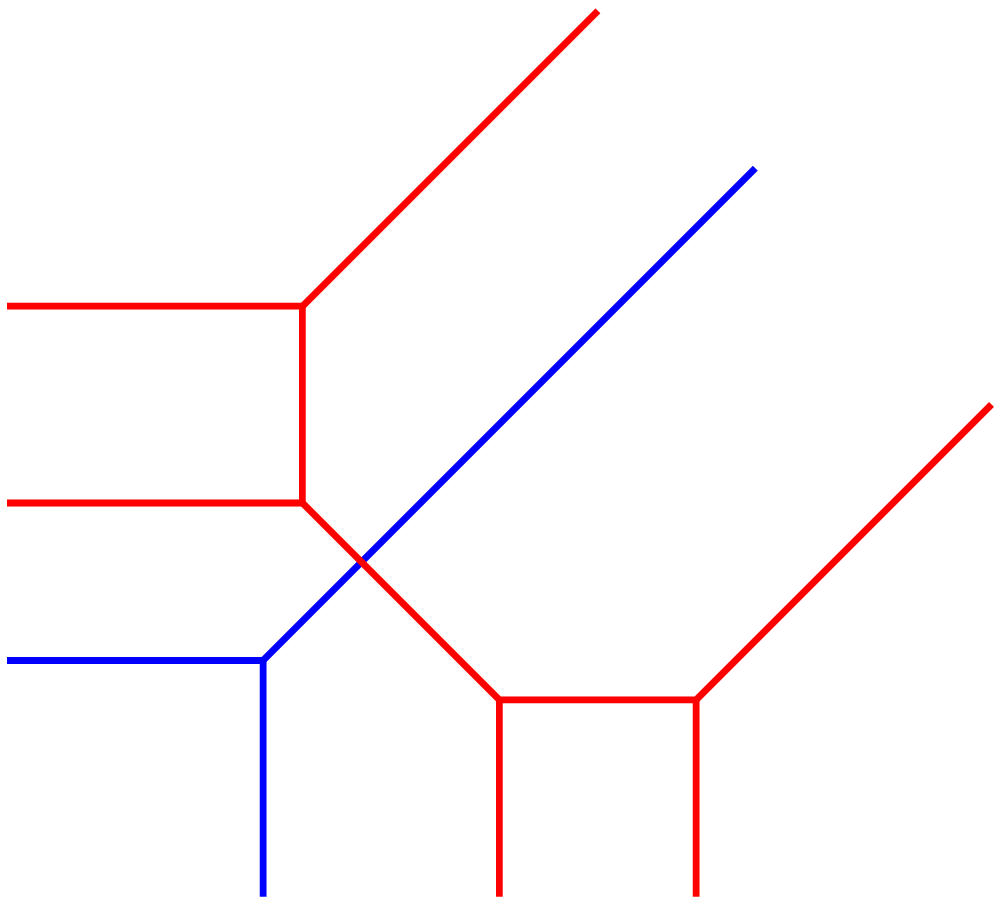} 
\\ \\a) && b) &&c)

\end{tabular}
\end{center}
\caption{Intersections de droites et de coniques tropicales}
\label{inter}
\end{figure}

Sauf accident, deux droites tropicales se coupent en un seul point
(voir  figure \ref{inter}a), comme en géométrie classique.
Maintenant, une droite et une conique tropicales se coupent-elles en
deux points?  
Si on compte naïvement le nombre de points d'intersection, la réponse
est: des fois oui (figure \ref{inter}b) des fois non (figure
\ref{inter}c)... 

En fait, l'unique point d'intersection de la conique et de la droite tropicales
sur la figure 
\ref{inter}c, doit être compté 2 fois. Mais pourquoi 2
ici et 1 dans les  cas précédents? La réponse se trouve dans
la subdivision duale à l'union des deux courbes. 

Remarquons tout d'abord que l'union de deux courbes tropicales $C_1$
et $C_2$ est 
encore une courbe tropicale. En effet, on vérifie aisément que l'union
de deux graphes équilibrés est encore un graphe équilibré, mais on peut
aussi voir que si les courbes tropicales $C_1$ et $C_2$ sont
respectivement   définies par les
polynômes tropicaux $P_1(x,y)$ et $P_2(x,y)$, alors le polynôme $Q(x,y)=\tg
P_1(x,y)P_2(x,y) \td$ définit précisément la courbe $C_1\cup C_2$. De
plus, le degré de $C_1\cup C_2$ est la somme des  degrés de $C_1$ et
de $C_2$. Ainsi, cela a bien un sens de parler de la subdivision duale
à la courbe $C_1\cup C_2$.  

\begin{figure}[h]
\begin{center}
\begin{tabular}{ccccc}
\includegraphics[width=2cm, angle=0]{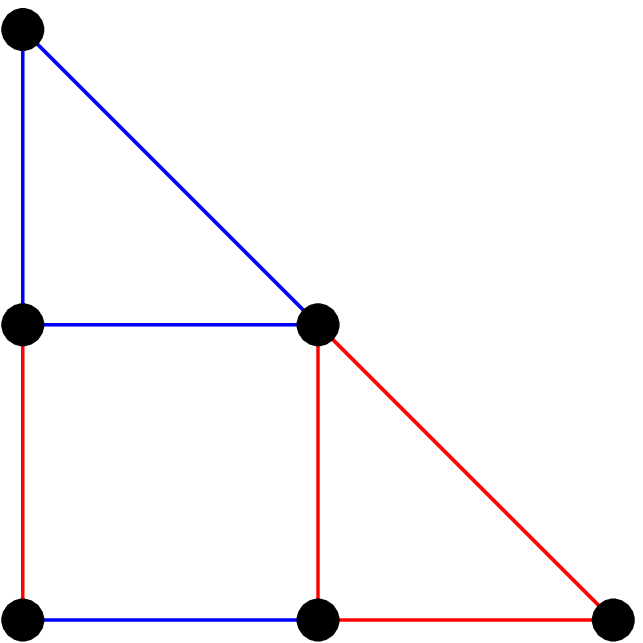}&\hspace{3ex} &
\includegraphics[width=3cm, angle=0]{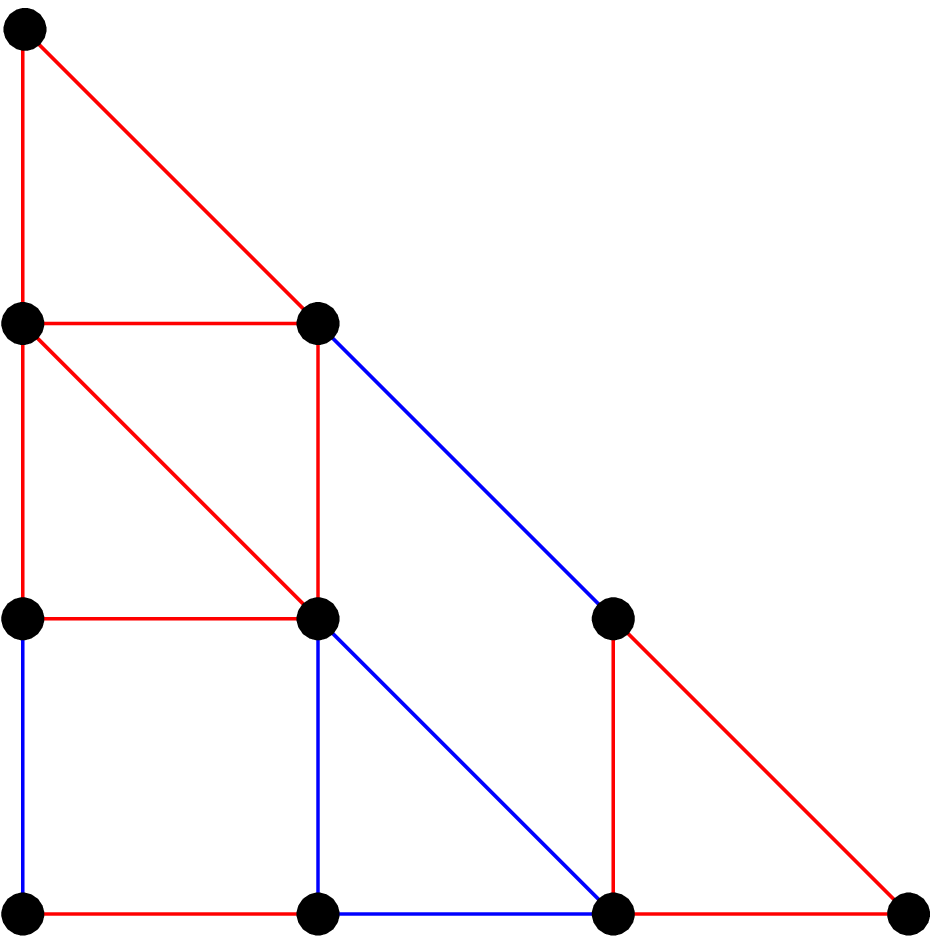}&\hspace{3ex} &
\includegraphics[width=3cm, angle=0]{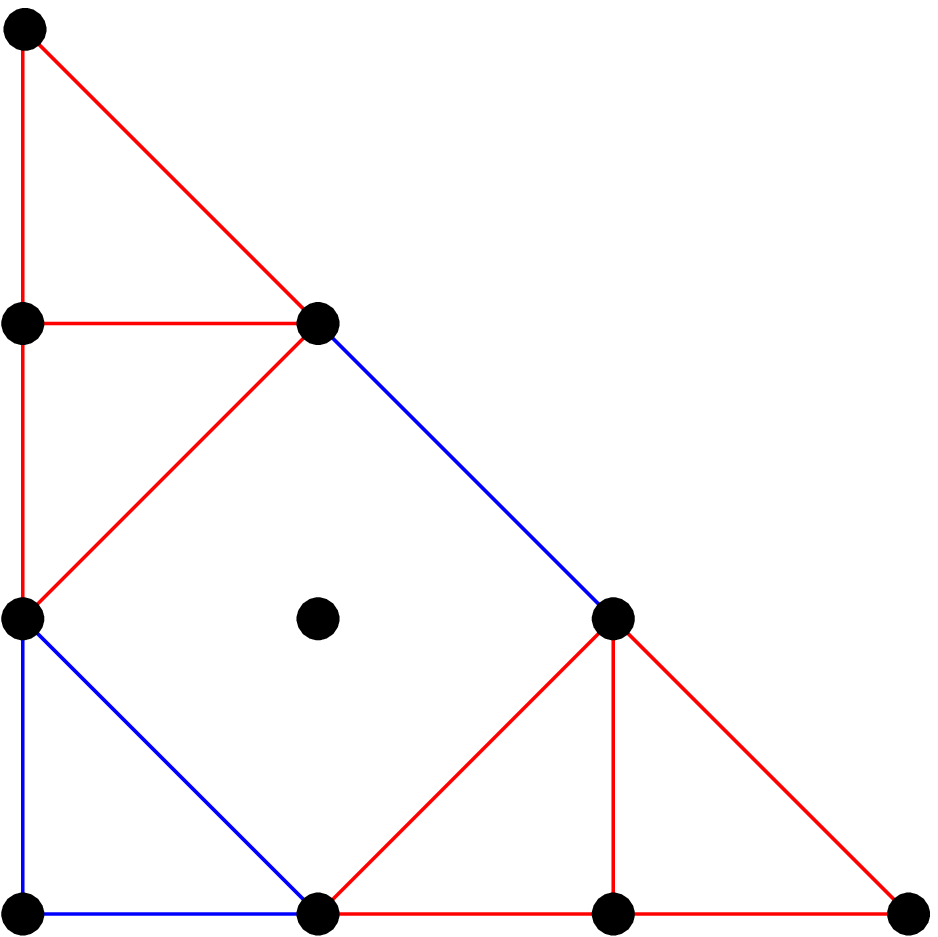} 
\\ \\a) && b) &&c)

\end{tabular}
\end{center}
\caption{Subdivisions duales à l'union des courbes de la figure \ref{inter}}
\label{sub inter}
\end{figure}

Les subdivisions duales à l'union des courbes $C_1$ et $C_2$ dans
chacun des cas de 
la figure \ref{inter} sont représentées sur la figure \ref{sub
  inter}. À chaque fois les sommets de $C_1\cup C_2$ 
sont les sommets de $C_1$, les sommets de $C_2$, et les points
d'intersection de $C_1$ et $C_2$. De plus comme chaque point de $C_1\cap
C_2$ est l'intersection d'une arête de $C_1$ et d'une arête de $C_2$,
le polygone dual à un tel sommet de $C_1\cup C_2$ est un parallélogramme.
Pour rendre la figure \ref{sub inter} plus transparente, nous avons dessiné
chaque arête de la subdivision duale de la même
couleur que son arête duale. 
Nous constatons alors que sur les
figures \ref{sub inter}a et b, les parallélogrammes correspondant sont
d'aire 1, alors que le parallélogramme de la subdivision de la figure
\ref{sub inter}c est d'aire 2! Ainsi, il semble que l'on doive compter
chaque point d'intersection avec la multiplicité définie ci-dessous.

\begin{defi}
Soit $C_1$ et $C_2$ deux courbes tropicales s'intersectant en un
nombre fini de points et en dehors des sommets 
des deux courbes, et soit $p$ un point
d'intersection de $C_1$ et $C_2$. La multiplicité
tropicale de $p$
comme point d'intersection de $C_1$ et $C_2$ est l'aire du
parallélogramme dual à $p$ dans la subdivision duale à $C_1\cup C_2$.
\end{defi}

Avec cette définition, démontrer le théorème de Bézout tropical est un
jeu d'enfant. 

\begin{thm}
Soit $C_1$ et $C_2$ deux courbes tropicales de degré $d_1$ et $d_2$
s'intersectant en un nombre fini de points et en dehors des sommets
des deux courbes. Alors la somme des multiplicités tropicales des points
d'intersection de $C_1$ et $C_2$ est égale à $d_1d_2$.
\end{thm}
\begin{proof}
Notons $s$ cette somme.
Il existe trois types de polygones dans la subdivision duale à la
courbe tropicale $C_1\cup C_2$: 
\begin{itemize}
\item ceux duaux à un sommet de $C_1$. La somme de leur aire est égale
  à l'aire de $\Delta_{d_1}$, c'est-à-dire $\frac{d_1^2}{2}$.
\item ceux duaux à un sommet de $C_2$. La somme de leur aire est égale
  à $\frac{d_2^2}{2}$.
\item ceux duaux à un point d'intersection de $C_1$ et $C_2$. La somme
  de leur aire est égale à $s$.
\end{itemize}
Puisque la courbe $C_1\cup C_2$ est de degré $d_1+d_2$, la somme de
l'aire de tous ces polygones est égale à l'aire de 
$\Delta_{d_1+d_2}$, c'est-à-dire $\frac{(d_1+d_2)^2}{2}$ et donc
$$s=\frac{(d_1+d_2)^2 -d_1^2 -d_2^2}{2}= d_1d_2  $$
\end{proof}

\subsection{Intersection stable}
\begin{figure}[h]
\begin{center}
\begin{tabular}{ccccccc}
\includegraphics[width=3cm, angle=0]{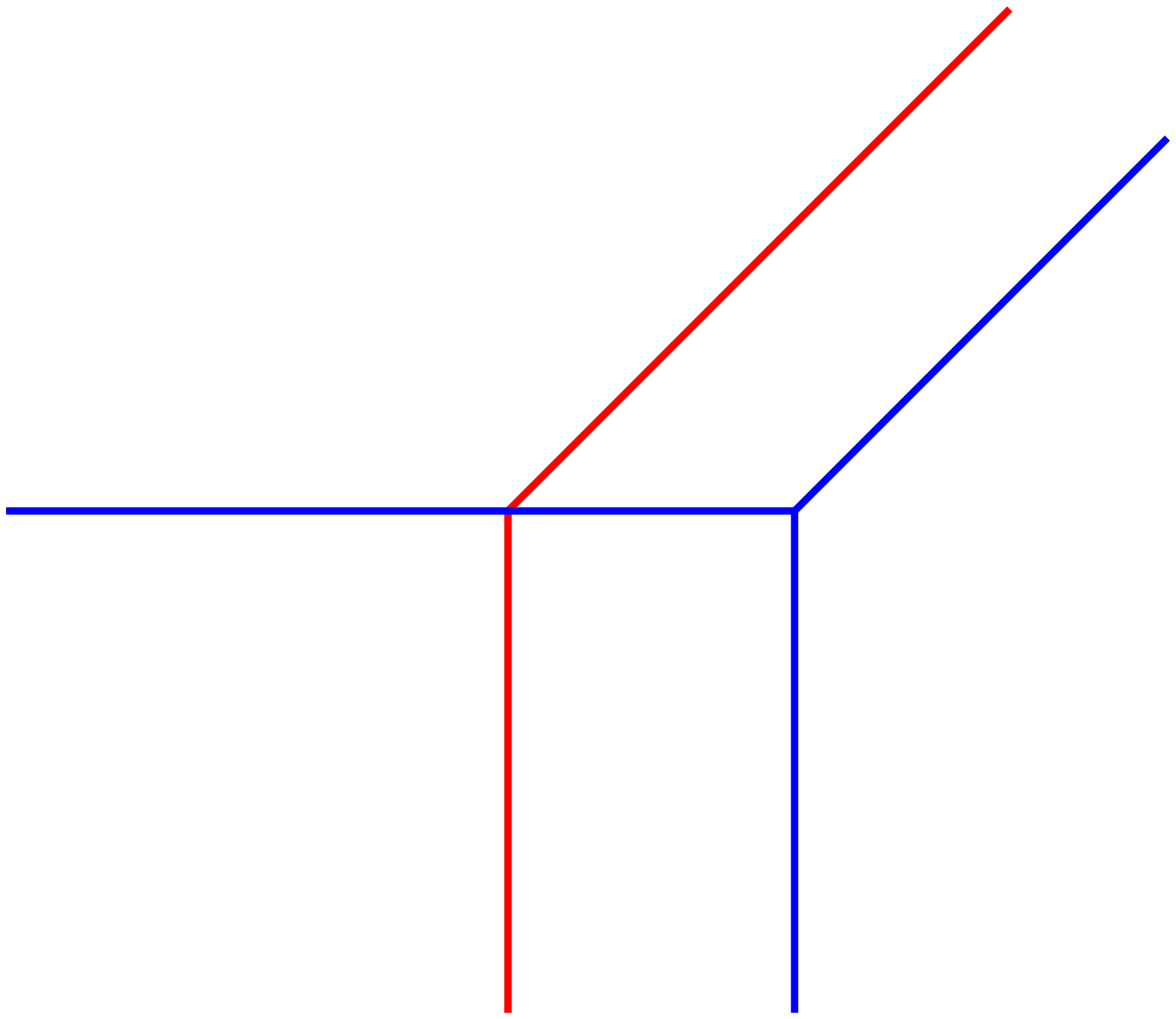}&\hspace{3ex} &
\includegraphics[width=3cm, angle=0]{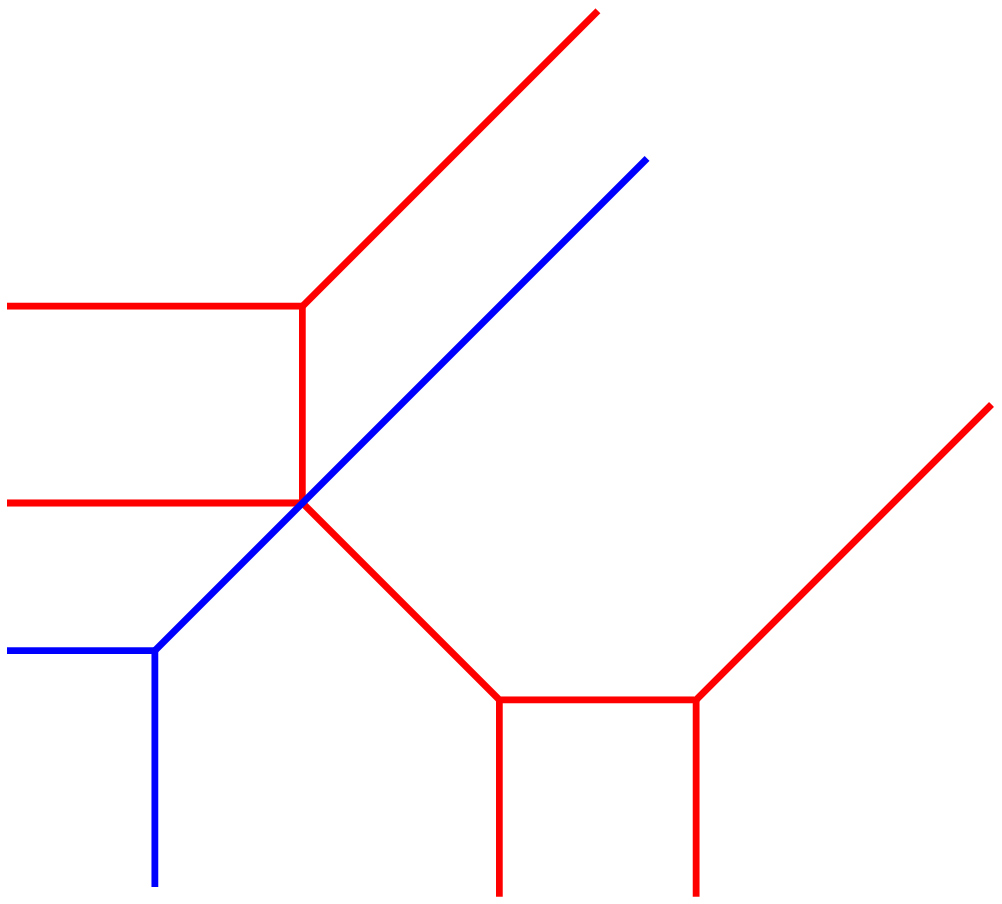}&\hspace{3ex} &
\includegraphics[width=3cm, angle=0]{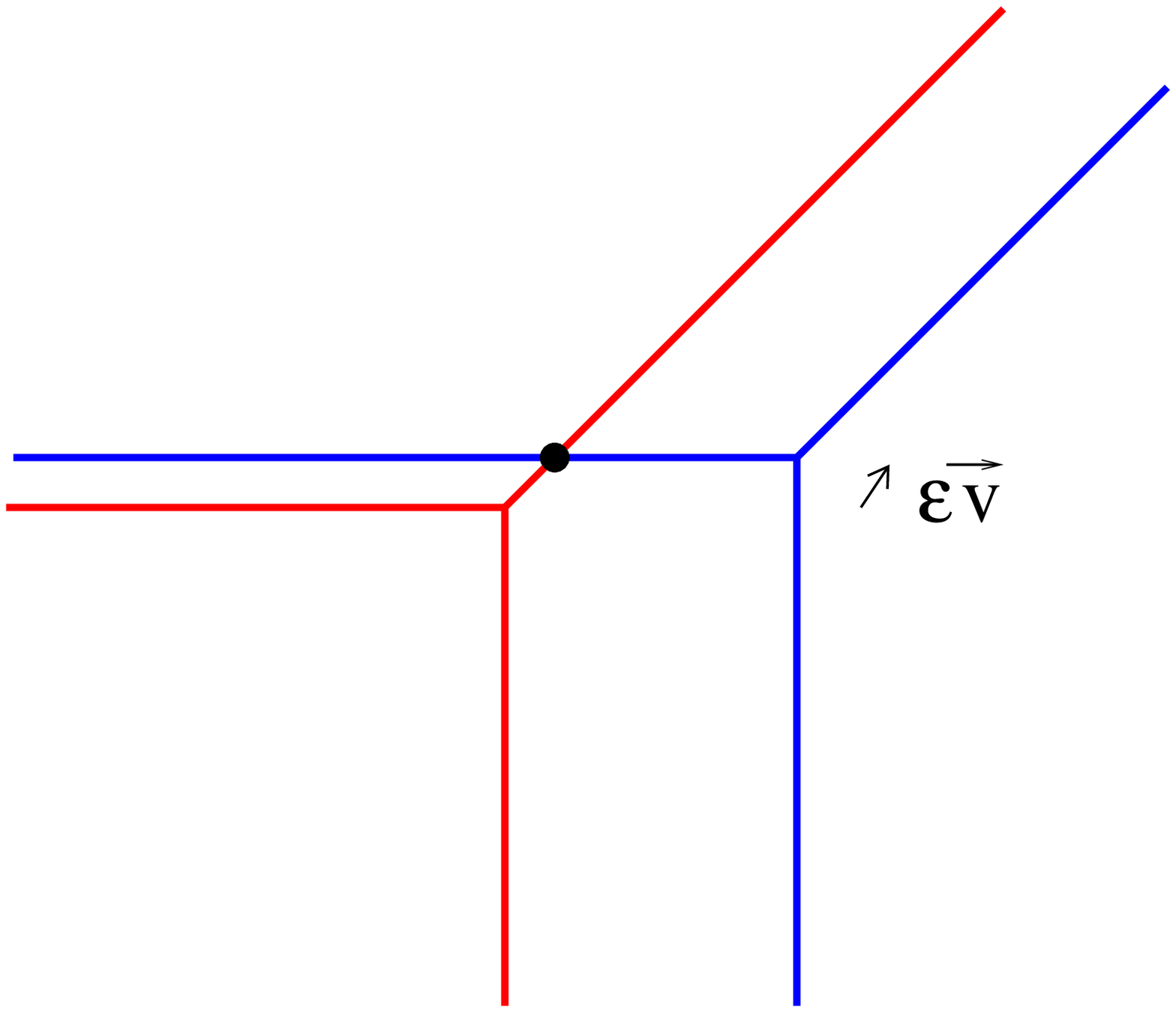}&\hspace{3ex} &
\includegraphics[width=3cm, angle=0]{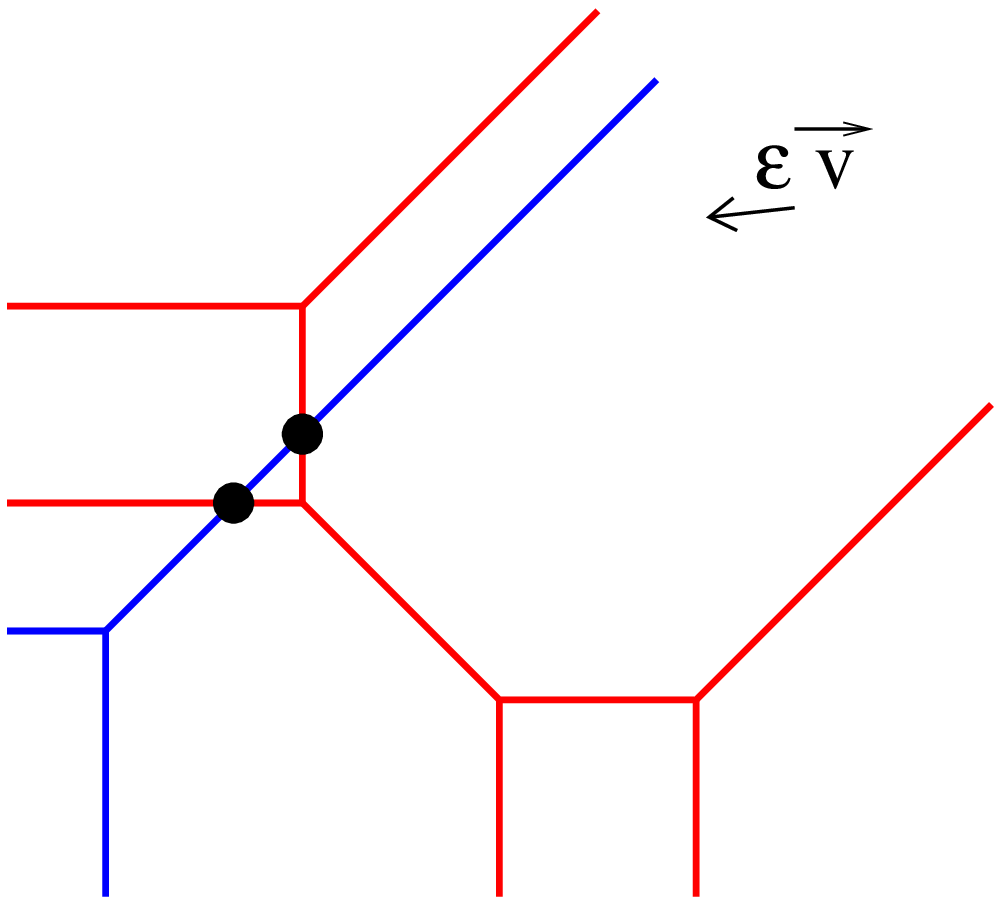} 
\\ \\a) && b)  &&c) &&d)

\end{tabular}
\end{center}
\caption{Intersections non-transverses et translation}
\label{inter st}
\end{figure}

Nous avons considéré jusqu'ici uniquement des courbes tropicales
s'intersectant ``gentiment'', c'est-à-dire en un nombre fini de points
et en dehors
des sommets. Mais que peut-on dire dans les deux cas représentés
sur les figures \ref{inter st}a (deux droites tropicales
s'intersectant le long d'une arête) et b (une droite passant par le
sommet d'une conique)? Heureusement, nous avons plus d'un tour tropical dans notre
poche. 

Soit $\varepsilon$  un petit nombre réel et $\vec v$  un
vecteur  dont le rapport des deux coordonnées
est un nombre irrationnel.
Si on translate dans chaque cas une des deux courbes du vecteur $\varepsilon
\vec v$, on se retrouve alors dans le cas d'une intersection gentille 
 (voir figures  \ref{inter st}c et d). Évidemment nos nouveaux
points d'intersection dépendent du vecteur $\varepsilon
\vec v$. En revanche  la limite de ces points lorsque
$\varepsilon$ tend vers 0 ne dépend pas de $\vec v$, ce sont les
\textit{points d'intersection stables} des deux courbes. Leur
multiplicité est égale à la somme des multiplicités des points
d'intersection dont
ils sont limite. 

Par exemple, le point d'intersection stable des deux droites de la figure
\ref{inter st}a est  le sommet de la droite de gauche, de
multiplicité 1. Nos deux droites 
tropicales se coupent  bien en un point.
Le point d'intersection stable des deux courbes de la figure
\ref{inter st}b est le sommet de la conique, de multiplicité 2.

Notons que l'intersection stable de deux courbes tropicales est
concentrée aux points d'intersection isolés et aux
sommets des deux courbes.  
Grâce à  l'intersection stable, nous pouvons supprimer du Théorème de
Bézout tropical les hypothèses de bonne position des deux courbes.

\begin{thm}
Soit $C_1$ et $C_2$ deux courbes tropicales de degré $d_1$ et
$d_2$. Alors la somme des multiplicités  des points 
d'intersection stables de $C_1$ et $C_2$ est égale à $d_1d_2$.
\end{thm}

Nous découvrons au passage un
phénomène tropical surprenant: une courbe tropicale a une
\textit{auto-intersection} bien définie\footnote{En géométrie
  algébrique classique, seul le nombre de points d'auto-intersection
  d'une courbe plane est 
  défini, pas leur position sur la
  courbe. Une droite s'auto-intersecte en 1 point, mais ce n'est pas
  clair quel est ce point...}! En effet il suffit simplement de considérer
 l'intersection 
stable de cette courbe tropicale avec elle même.
D'après ce qui précède, cette auto-intersection est concentrée aux
sommets de la courbe (voir figure \ref{auto
  inter}).

\begin{figure}[h]
\begin{center}
\begin{tabular}{c}
\includegraphics[width=10cm, angle=0]{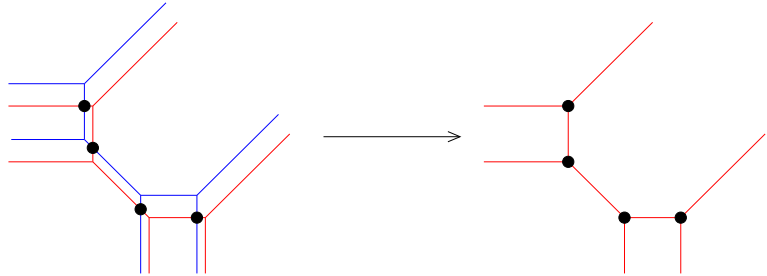}

\end{tabular}
\end{center}
\caption{4 points d'auto-intersection d'une conique tropicale}
\label{auto inter}
\end{figure}

\subsection{Exercices}
\begin{exo}
\begin{enumerate}
\item Déterminer les points d'intersection stables des deux courbes
  tropicales de l'exercice 1 de la partie 2, ainsi que leur multiplicité.
\item Un point double d'une courbe tropicale est l'intersection de
  deux arêtes de celle-ci. Montrer qu'une conique tropicale avec un
  point double est l'union de deux droites 
  tropicales. On pourra considérer une droite passant par
    le point double et un autre sommet de la conique.
\item Montrer qu'une courbe tropicale de degré 3 avec deux points
  doubles est l'union d'une droite et d'une conique tropicale. Montrer
  qu'une courbe tropicale de degré 3 avec trois points 
  doubles est l'union de trois droites 
  tropicales. 
\end{enumerate}
\end{exo}

\section{Quelques explications}\label{expl}

Arrêtons nous quelques instants dans l'étude 
de la géométrie tropicale
proprement dite, et donnons quelques raisons du lien très fort entre
géométrie classique et géométrie tropicale. Notre but est d'illustrer en
particulier le fait que la géométrie tropicale est une limite de la
géométrie classique.
 Pour résumer grossièrement le contenu
de cette partie,  la géométrie
tropicale est l'image de la géométrie classique par le logarithme en
base $\infty$.

\subsection{Déquantification de Maslov}
Expliquons d'abord comment le semi-corps tropical arrive naturellement
comme limite de semi-corps classiques. Ce 
processus, étudié par Victor Maslov et
ses collaborateurs à partir des années 90, s'appelle la \textit{déquantification
  des nombres réels}.

Un semi-corps bien connu est $(\RR_+,+,\times)$, l'ensemble des nombres
réels
positifs ou nuls muni de l'addition et de la multiplication
classiques.
Si $t$ est un nombre strictement positif, alors le logarithme en base
$t$ fournit une bijection entre $\RR_+$ et $\TT$, et cette bijection
induit une structure de semi-corps sur $\TT$ où les opérations, notées
$\tg +_t\td$ et $\tg\times_t\td$, sont données par:
$$\tg x +_t  y\td= \log_t(t^x +t^y) \  \ \ \text{ et } \ \ \  \tg x \times_t  y\td=
\log_t(t^x t^y) = x+y$$ 

On voit donc déjà apparaître l'addition classique comme une multiplication
exotique sur $\TT$.
Notons que par construction, tous les semi-corps $(\TT,\tg +_t\td,\tg
\times_t\td)$ sont isomorphes à $(\RR_+,+,\times)$. L'inégalité
triviale $\max(x,y)\le x+y\le 2\max(x,y)$ sur $\RR_+$ combinée à la
croissance de la fonction logarithme nous donne l'encadrement suivant:
$$\forall t>0,  \ \  \max(x,y) \le \tg x \times_t  y\td \le \max(x,y)
+\log_t 2$$
Lorsque $t$ tend vers l'infini $\log_t 2$ tend vers $0$, et la loi
$\tg +_t\td$ tend donc vers 
l'addition tropicale $\tg +\td$! Ainsi, le semi-corps tropical arrive
naturellement 
comme 
 dégénérescence du semi-corps classique
$(\RR_+,+,\times)$. Ou encore, on peut voir le semi-corps classique
$(\RR_+,+,\times)$ comme une déformation du semi-corps tropical, d'où
 le terme de ``déquantification''.

\subsection{Déquantification d'une droite du plan}\label{deqdte}
\begin{figure}[h]
\begin{center}
\begin{tabular}{cccccc}
\includegraphics[width=2.5cm,  angle=0]{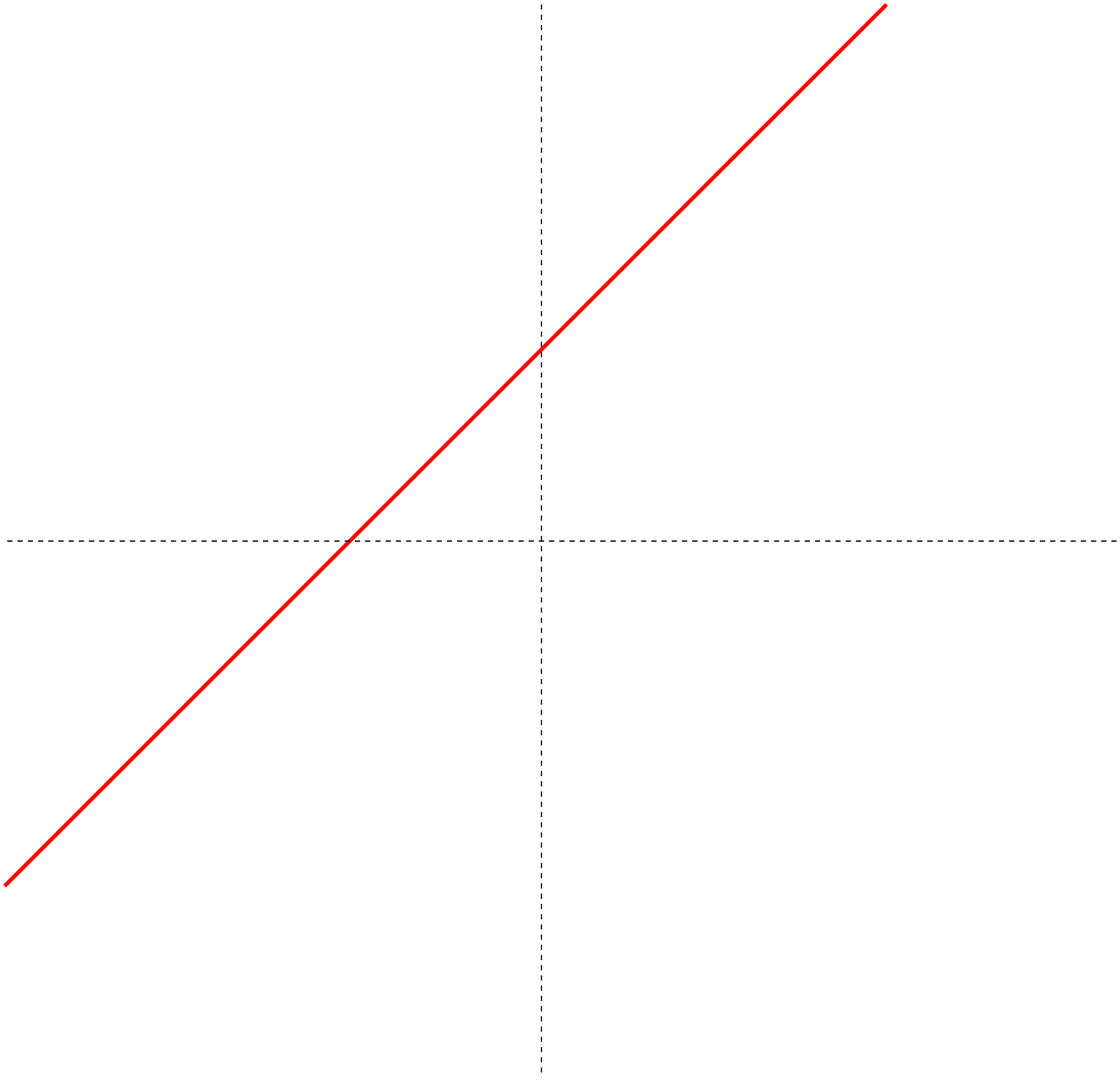}&
\includegraphics[width=2.5cm, angle=0]{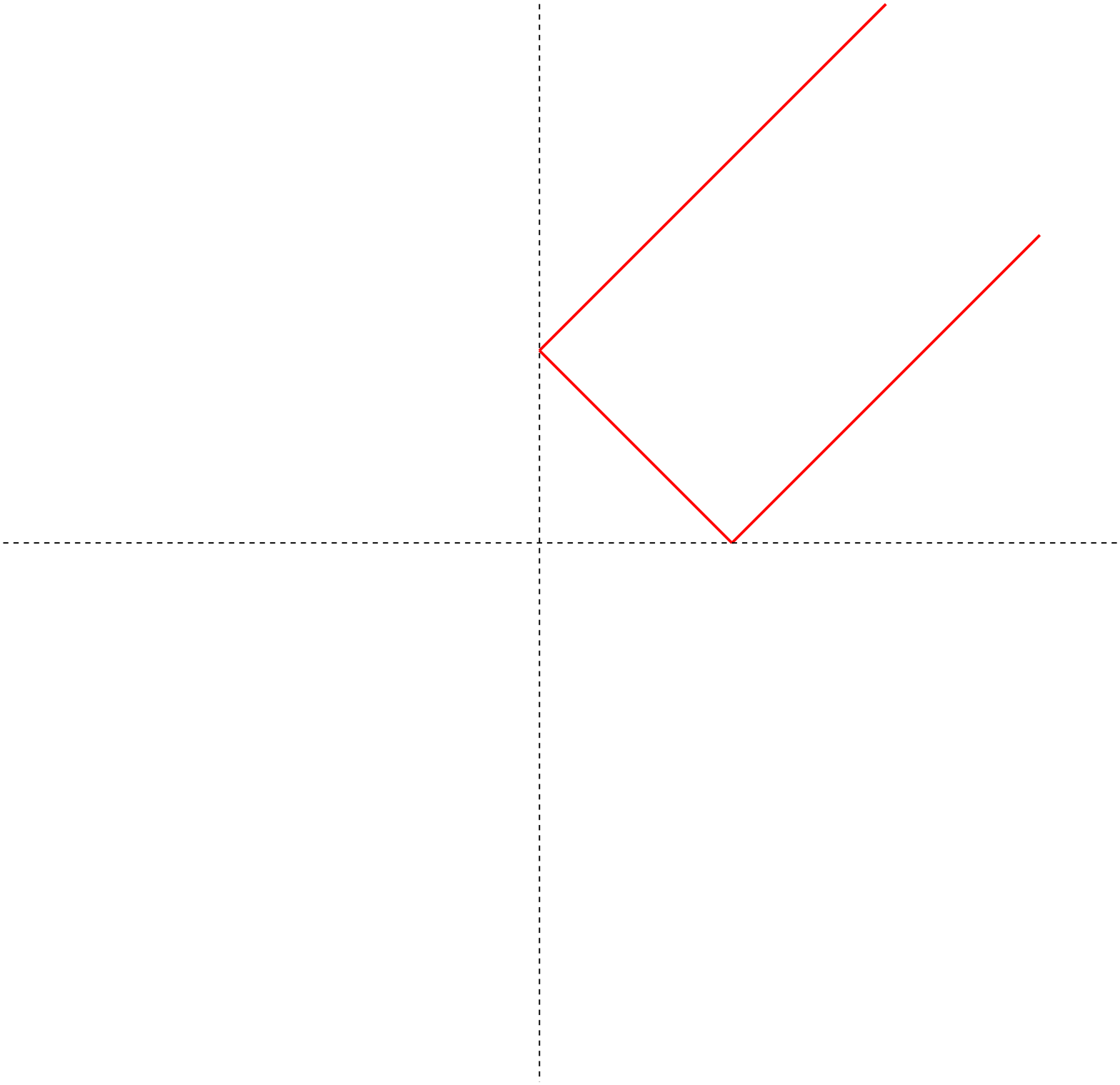}&
\includegraphics[width=2.5cm, angle=0]{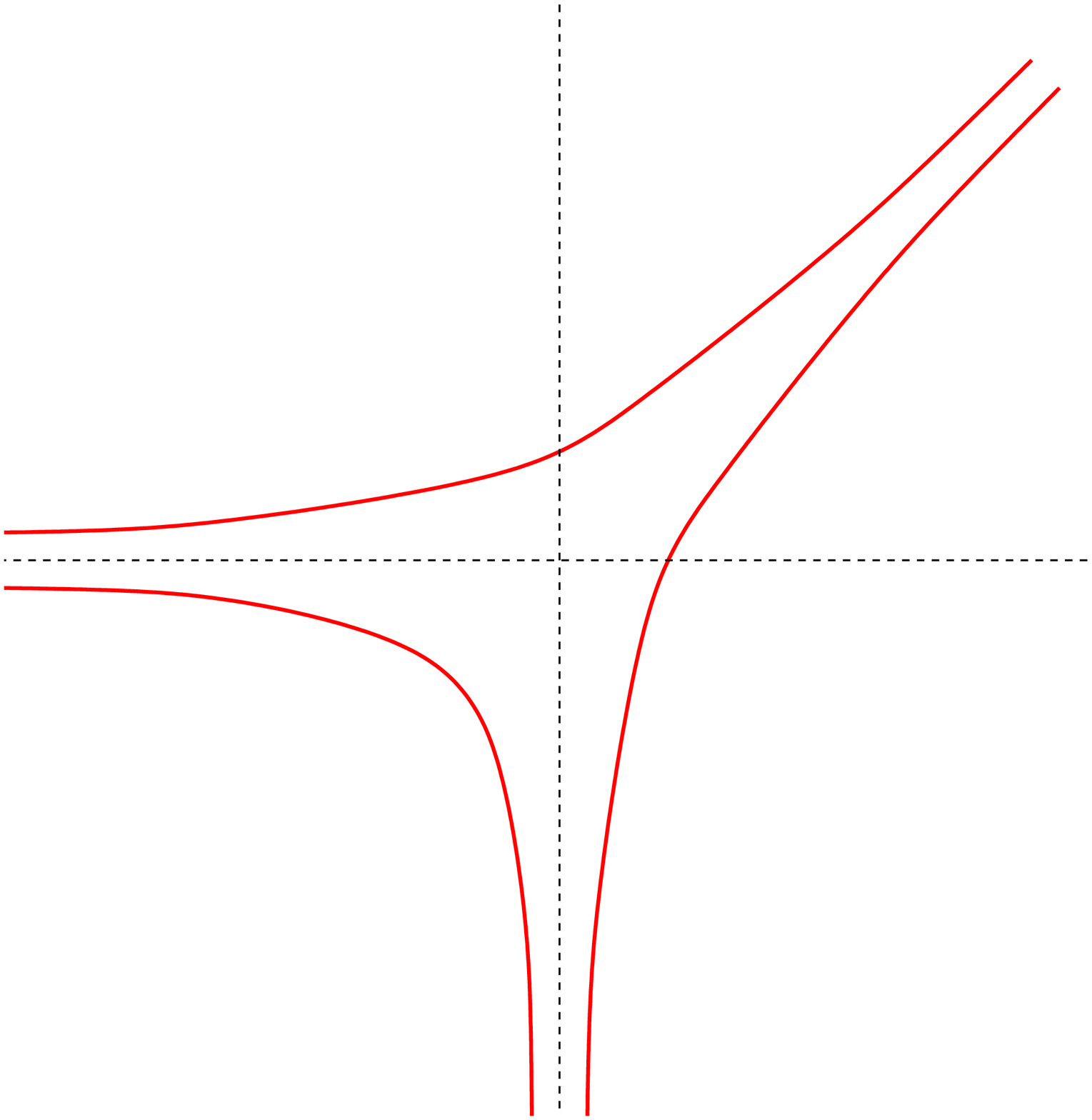}&
\includegraphics[width=2.5cm, angle=0]{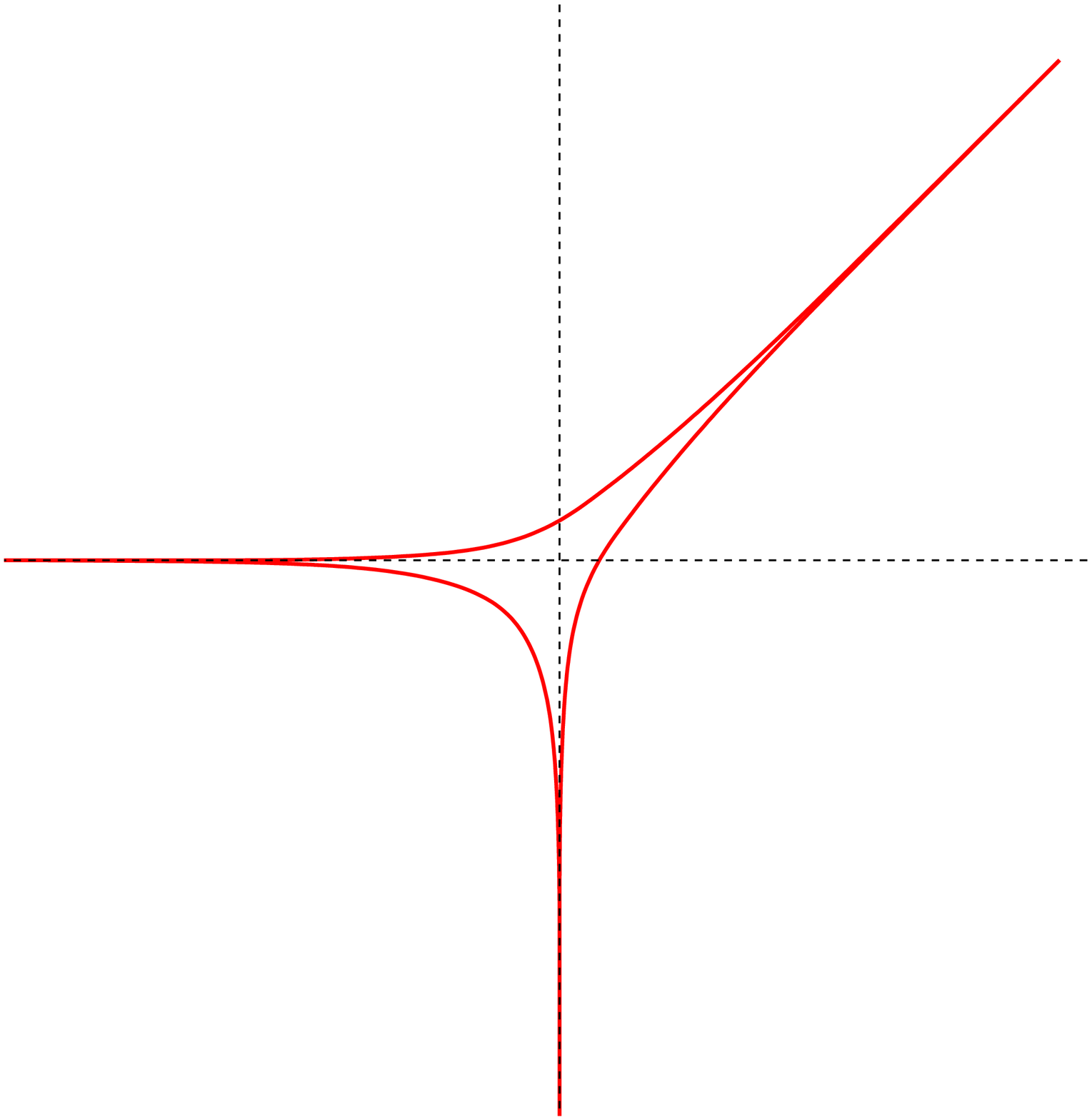}&
\includegraphics[width=2.5cm, angle=0]{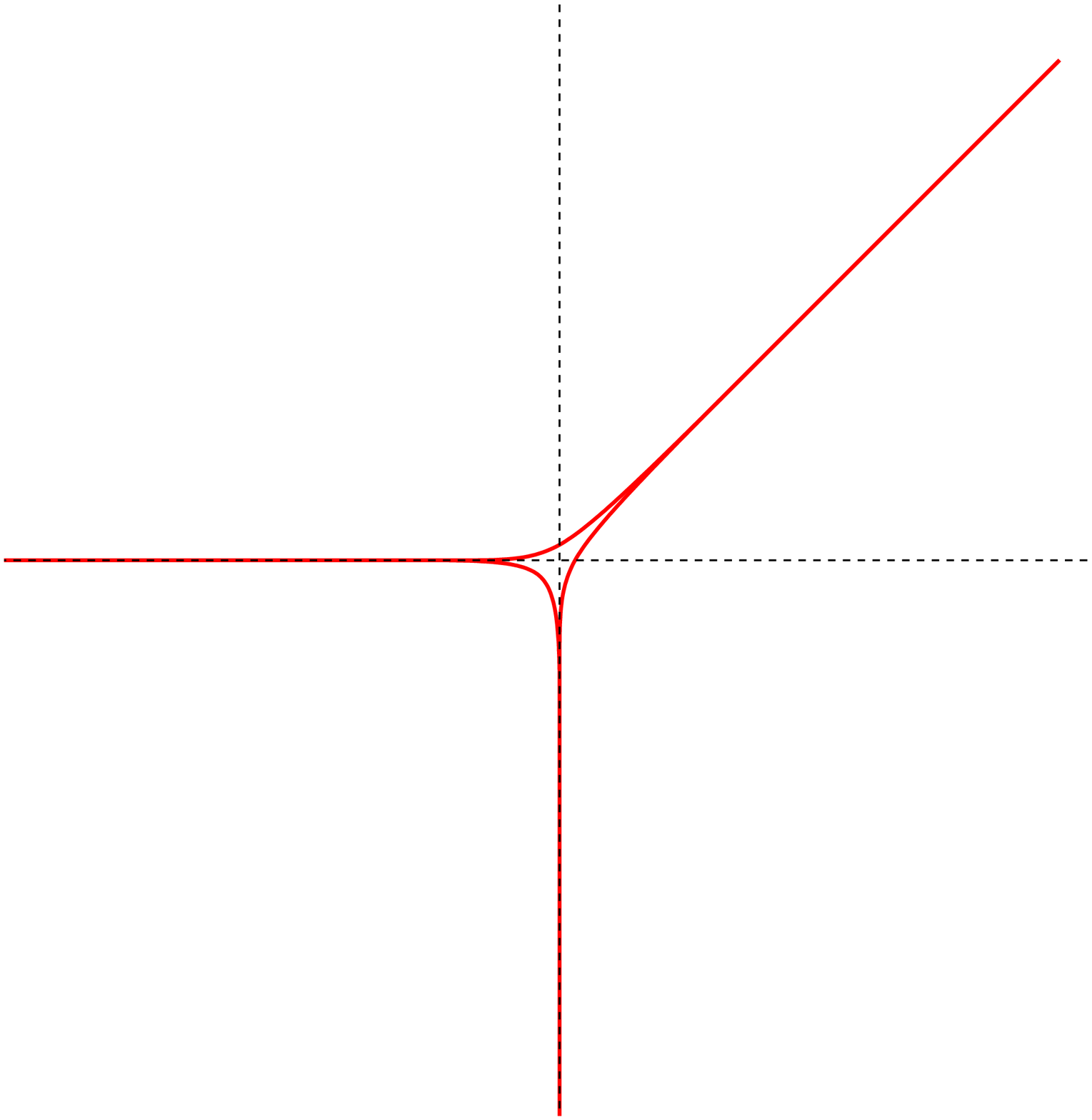}&
\includegraphics[width=2.5cm, angle=0]{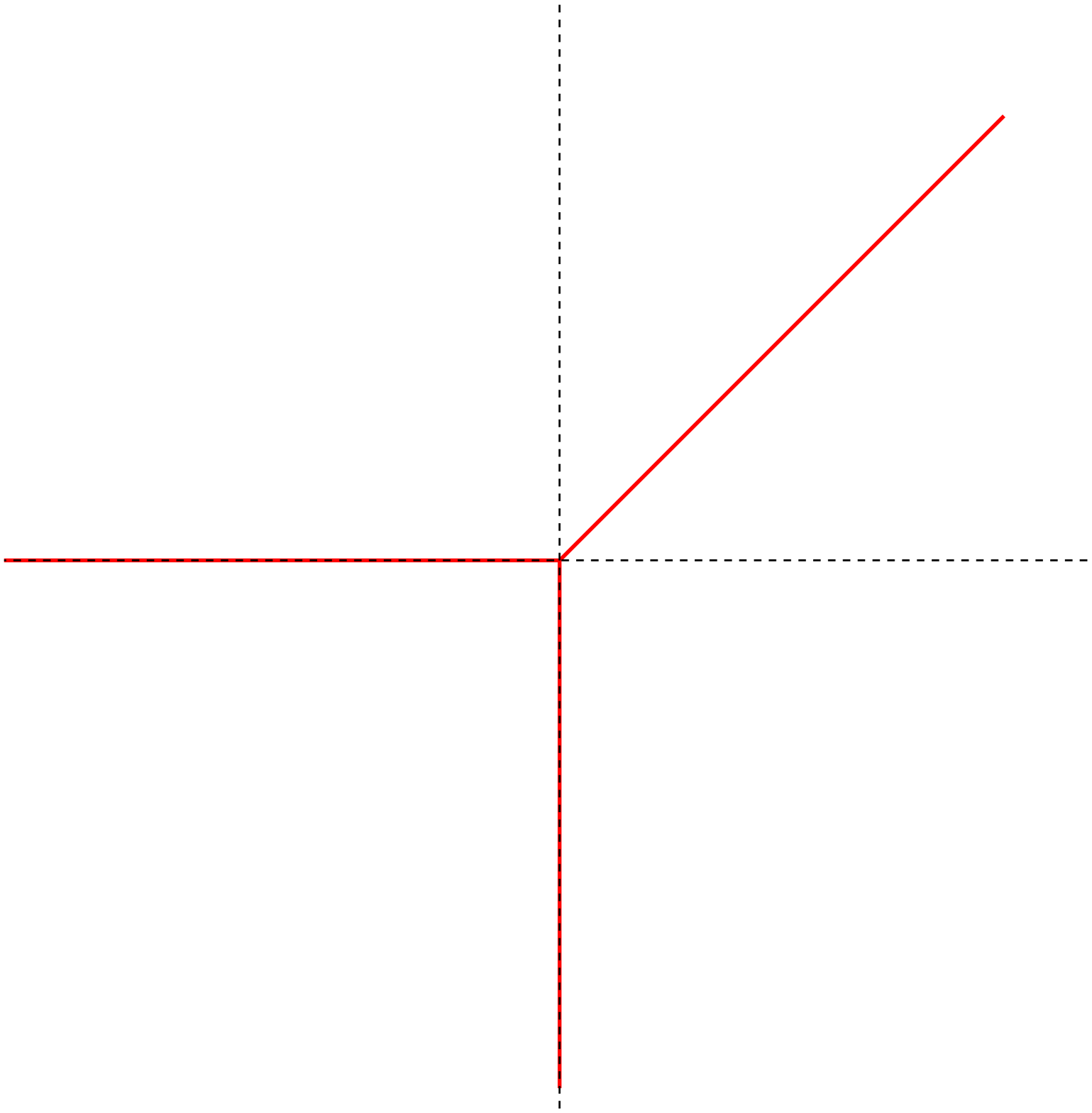}
\\ \\ a) & b) & c) &d)&e)& f)

\end{tabular}
\end{center}
\caption{Déquantification d'une droite}
\label{amibe}
\end{figure}

Appliquons un raisonnement similaire avec la droite d'équation $x-y+1$
dans le plan $\RR^2$ (voir  figure \ref{amibe}a). 
Replions tout d'abord les 4 quadrants sur le quadrant positif par
l'application valeur absolue (voir  figure
\ref{amibe}b). L'image par le logarithme en base $t$ de cette droite
repliée dans
$(\RR_+^*)^2$ ressemble au dessin de la figure \ref{amibe}c. Par définition,
prendre le logarithme en base $t$ revient à prendre le logarithme
népérien puis à appliquer une homothétie de rapport
$\frac{1}{\ln t}$. Ainsi, lorsque $t$ augmente,  l'image par le logarithme
en base $t$ de la valeur
absolue de notre droite se concentre sur un voisinage de l'origine et
des 3 directions asymptotiques (voir figures \ref{amibe}c, d et e). Et lorsqu'on fait tendre $t$ vers l'infini,
on voit apparaître sur la figure \ref{amibe}f... une droite tropicale!

\section{Patchwork}\label{patchwork}
La figure \ref{amibe} lue de gauche à droite nous explique comment
partir d'une droite classique dans le plan pour arriver à une droite
tropicale. 
La lecture de cette figure de droite à gauche est en fait beaucoup
plus intéressante! En effet,  nous voyons ainsi comment construire une
droite classique à partir d'une droite tropicale. 
La technique
dite du
\textit{patchwork} est une généralisation de cette observation.
 En particulier, elle fournit un  procédé purement
combinatoire
de construction de courbes algébriques
réelles à partir des courbes tropicales.
Mais avant
d'expliquer cette méthode en détail, remontons un  peu dans le passé.

\subsection{Le 16ème problème de Hilbert}

\begin{figure}[h]
\begin{center}
\begin{tabular}{ccccccc}
\includegraphics[width=3cm,  angle=0]{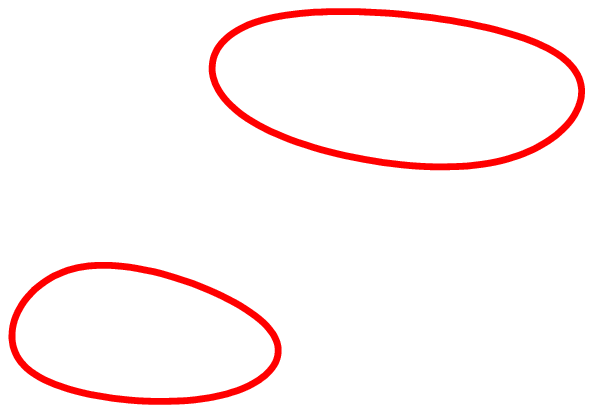}&\hspace{3ex} &
\includegraphics[width=3cm, angle=0]{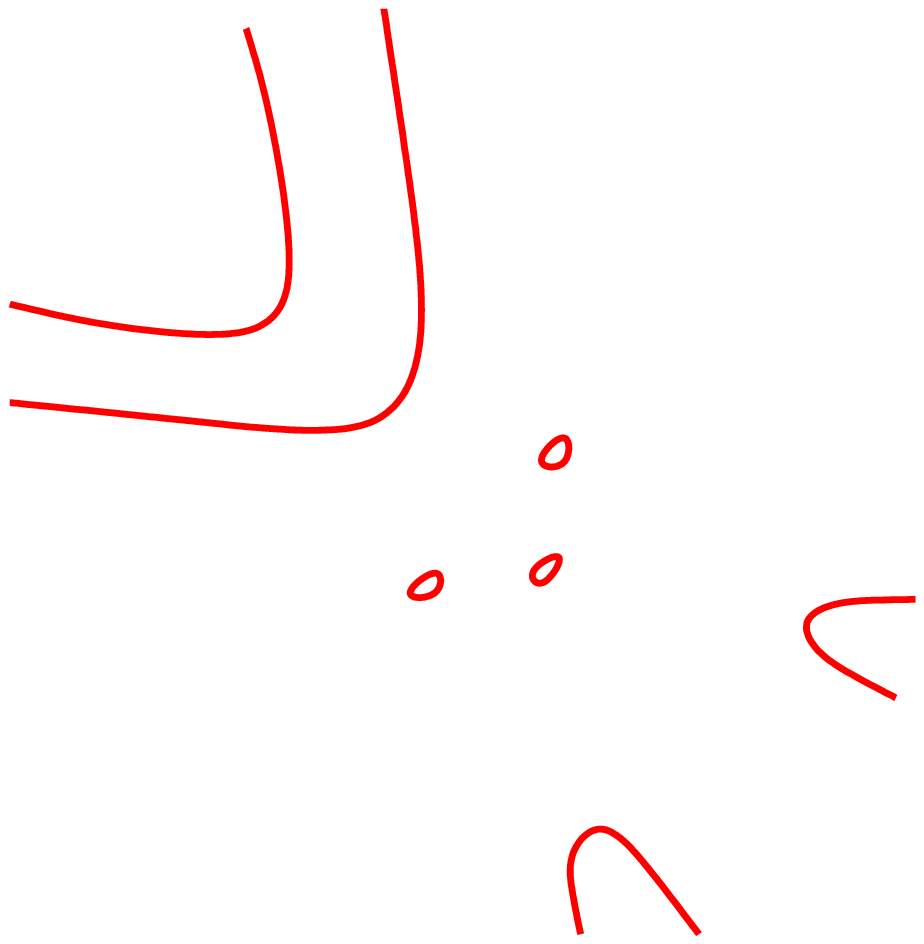}&\hspace{3ex} &
\includegraphics[width=3cm, angle=0]{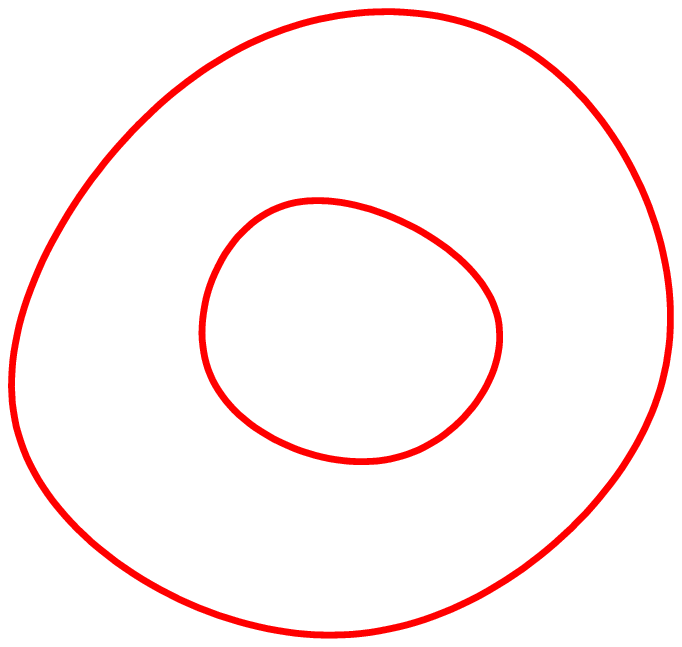}&\hspace{3ex} &
\includegraphics[width=3cm, angle=0]{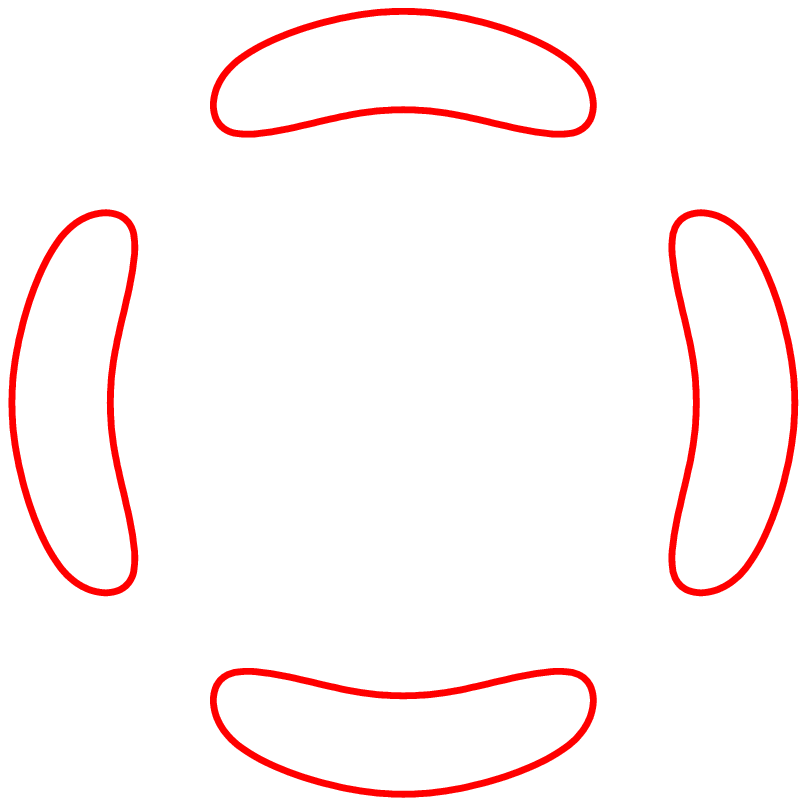}
\\ \\ a) && b) && c) &&d)

\end{tabular}
\end{center}
\caption{Quelques courbes algébriques réelles de degré 4}
\label{quartic}
\end{figure}

Une \textit{courbe algébrique réelle plane} est une courbe du plan
$\RR^2$ 
définie par une équation de la forme $P(x,y)=0$ où $P(x,y)$ est un polynôme à
coefficients réels. Les courbes algébriques réelles de degré 1 et 2
sont  simples et bien connues,
ce sont les droites et les coniques.
 À mesure que le degré de $P(x,y)$ augmente, le dessin réalisé par la courbe 
d'équation $P(x,y)=0$ peut être de plus en plus complexe.
 Pour s'en convaincre, il suffit de jeter un  \oe il à la figure
 \ref{quartic} où sont représentés une partie 
 des dessins possibles réalisés par une courbe algébrique réelle de degré 4.

Un théorème dû à Axel Harnack à la fin du XIXème siècle affirme qu'une
courbe algébrique réelle 
plane de degré $d$ a au maximum $\frac{d(d-1)+2}{2}$
composantes connexes. Mais comment  ces composantes 
peuvent-elles se situer les unes par rapport aux autres?
On appelle \textit{arrangement} d'une courbe algébrique
 réelle plane la position relative de ses composantes connexes dans le
 plan. Autrement dit, on ne s'intéresse pas à la position exacte de la
 courbe dans le plan, mais seulement au dessin qu'elle réalise. Par
 exemple si une courbe a deux composantes connexes bornées, on  se préoccupe
 uniquement de savoir si ces composantes
 sont en
 dehors l'une de l'autre (voir figure \ref{quartic}a) ou pas (voir
 figure \ref{quartic}c). 
Lors du deuxième congrès
international de mathématiques à Paris en 1900, David Hilbert énonça sa
célèbre liste de 23 problèmes pour le XXème siècle, et la première
partie de son 16ème problème peut être compris dans sa forme (très)
étendue comme 
suit:
\begin{center}
\textit{Étant donné un entier $d$, établir la liste des arrangements
  possibles  des courbes algébriques réelles de degré $d$.}
\end{center}

 À l'époque de Hilbert,  la
réponse était connue pour les courbes de degré au plus 4. 
Malgré des avancées spectaculaire dans ce problème au XXème siècle,
dues en particulier à 
des mathématiciens de l'école russe, de nombreuses questions restent
toujours sans réponses\footnote{Un problème plus raisonnable 
  et naturel consiste à regarder les arrangements des composantes
  connexes des courbes algébriques réelles projectives
  non-singulières. Pour ce problème plus restreint, la réponse est
  actuellement connue jusqu'en degré 7. C'est un théorème d'Oleg  Viro, et
  le patchwork est un outil essentiel de sa démonstration.}...

\subsection{Courbes réelles et courbes tropicales}

C'est en général un problème  difficile de construire une courbe
algébrique réelle d'un degré donné réalisant un arrangement
donné. Depuis plus d'un siècle, les mathématiciens ont proposé de
nombreuses et ingénieuses méthodes pour cela. 
Le patchwork inventée par
Oleg Viro dans les années 70  est une 
 des méthodes actuelle
les plus puissantes. À cette époque la géométrie tropicale
n'existait pas encore, et Viro énonça son théorème  dans un
langage différent du notre ici. Cependant il réalisa à la fin des
années 90 que le 
patchwork pouvait être interprété comme une \textit{quantification} des courbes
tropicales. Le patchwork, c'est donc
lire la figure \ref{amibe} de droite à gauche au lieu   de gauche à droite.
 Grâce à
cette interprétation nouvelle,
 Grigory Mikhalkin
généralisa peu après la méthode de Viro originelle.
Nous donnons ici une version simplifiée du
 patchwork, le lecteur intéressé trouvera une version plus complète
 dans les références indiquées à la partie \ref{ref}.

Dans la suite, si $a$ et $b$ sont deux nombres entiers nous notons
$s_{a,b}:\RR^2\to\RR^2$ la composée de $a$ symétries par rapport à l'axe des
abscisses suivies de $b$ symétries par rapport à l'axe de
ordonnées. Ainsi seules les valeurs modulo 2 de $a$ et $b$ sont
importantes, et $s_{0,0}$ est l'identité, $s_{1,0}$ est la
symétrie par rapport à l'axe des abscisses, $s_{0,1}$ est la
symétrie par rapport à l'axe des ordonnées et $s_{1,1}$ est la
symétrie par rapport à l'origine.

\begin{figure}[h]
\begin{center}
\begin{tabular}{ccccccc}
\includegraphics[width=3cm,  angle=0]{Figures/Droite.eps}&\hspace{3ex} &
\includegraphics[width=3cm, angle=0]{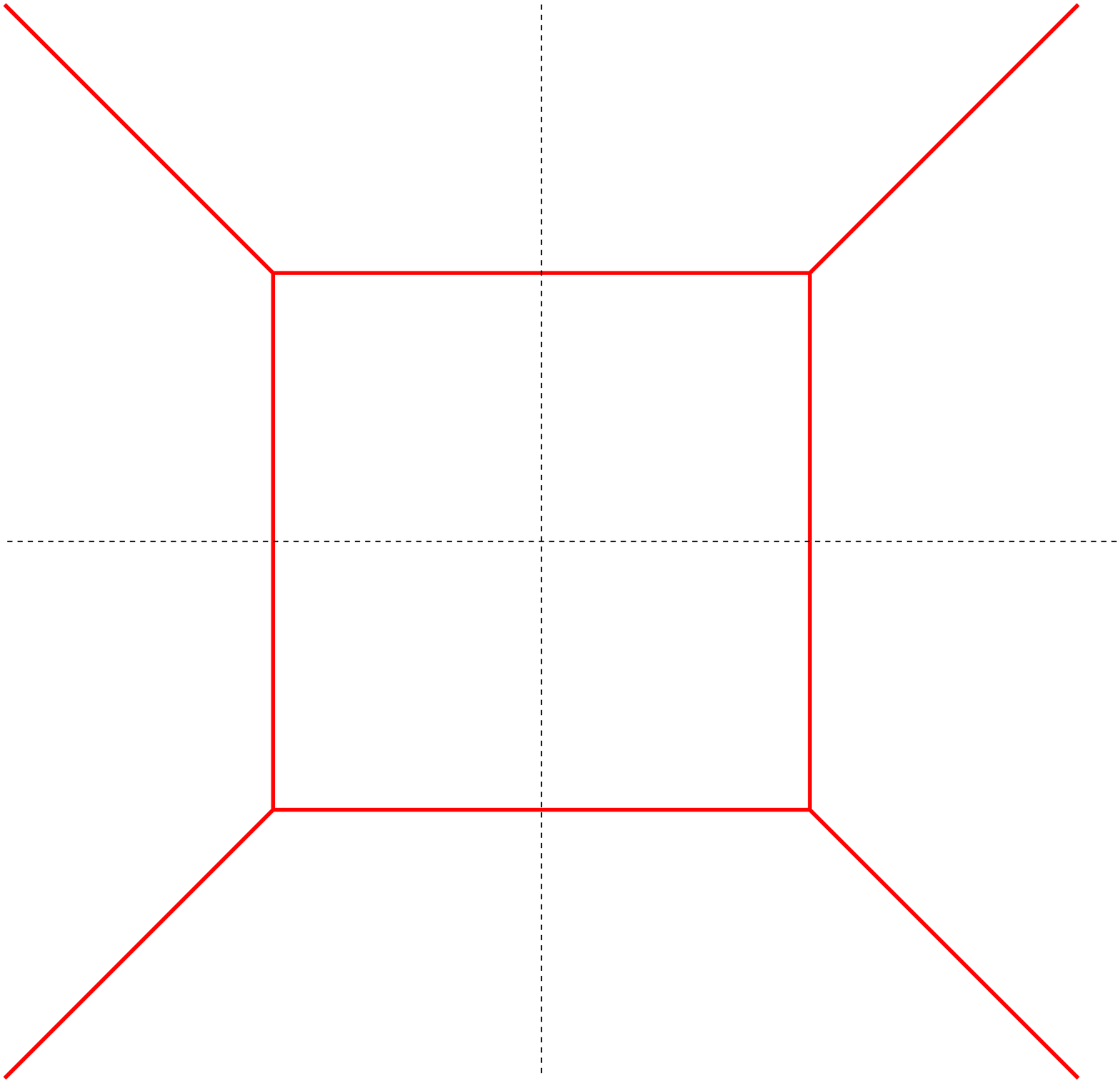}&\hspace{3ex} &
\includegraphics[width=3cm, angle=0]{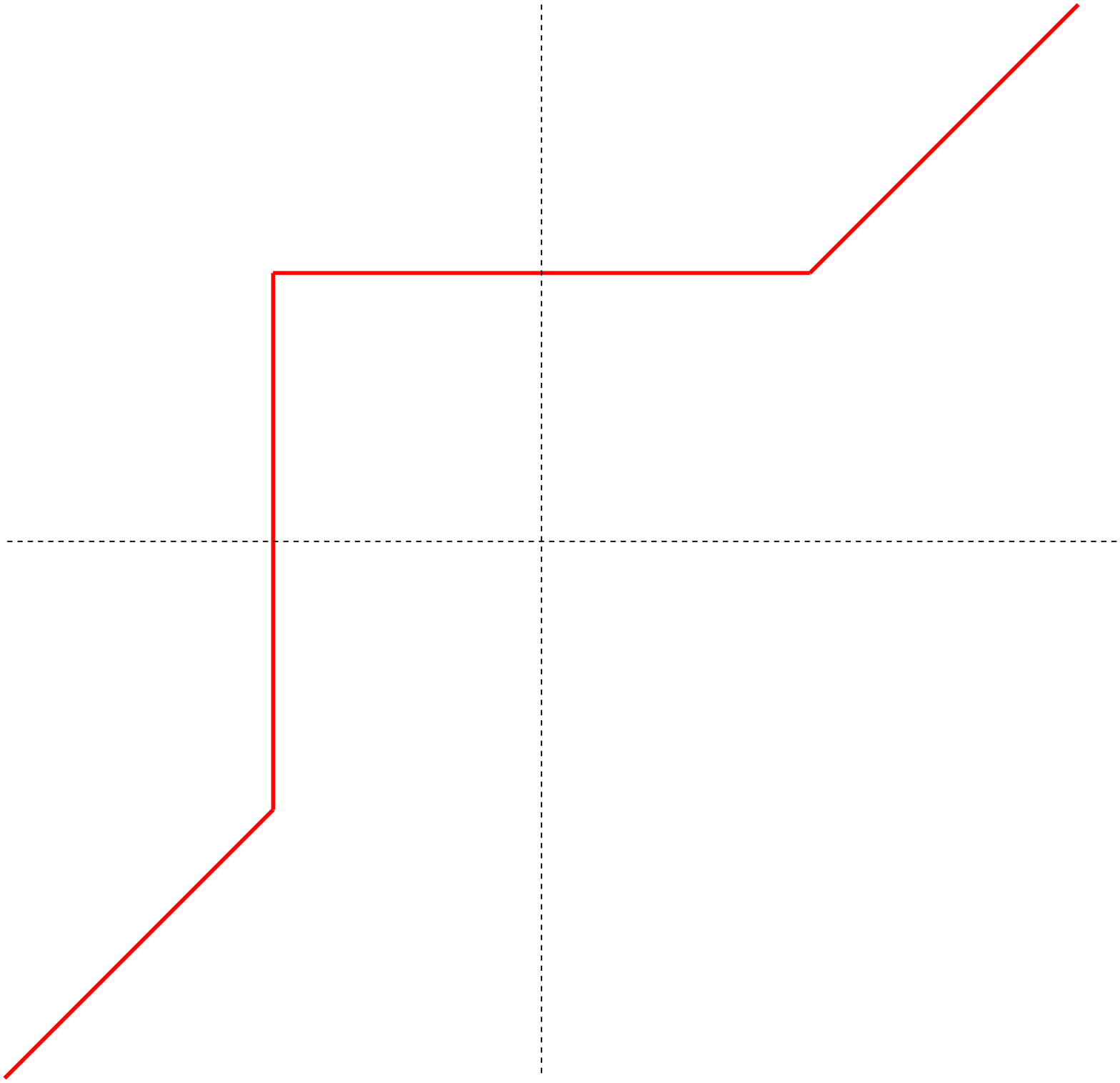}&\hspace{3ex} &
\includegraphics[width=3cm, angle=0]{Figures/PatchDte4.eps}
\\ \\ a) && b) && c) &&d)

\end{tabular}
\end{center}
\caption{Patchwork d'une droite}
\label{patch dte}
\end{figure}

Expliquons maintenant en détail la procédure de patchwork.
Prenons une courbe tropicale $C$ de degré $d$ n'ayant que
 des arêtes de poids impair et dont tous les
 polygones de la subdivision duale sont des triangles. Par exemple
 choisissons  la
 droite tropicale de la 
figure \ref{patch dte}a.  Pour chaque arête $e$ de $C$, choisissons un
vecteur $\vec v_e=(\alpha_e,\beta_e)$ directeur de $e$ avec $\alpha_e$
et $\beta_e$ deux 
entiers premiers entre eux. Pour la droite tropicale, choisissons les
vecteurs $(1,0)$, 
$(0,1)$ et $(1,1)$. 
Considérons maintenant que le $\RR^2$ dans
lequel vit notre courbe 
tropicale est en fait  le quadrant positif $(\RR^*_+)^2$ de $\RR^2$, et
prenons l'union de notre courbe tropicale avec ses copies
symétriques par rapport aux axes. Dans le cas de la droite tropicale,
nous obtenons alors la figure \ref{patch 
  dte}b.
Pour chaque arête $e$ de notre
courbe, effaçons  $e'$ et $e''$ deux des quatre
copies symétriques de $e$ suivant les deux règles suivantes:
\begin{itemize}
 \item $e'=s_{\alpha_e,\beta_e}(e'')$,
\item pour chaque sommet $v$ de $C$ adjacent aux arêtes $e_1$, $e_2$
  et $e_3$ et pour chaque couple
  $(\varepsilon_1,\varepsilon_2)$ dans $\{0,1\}^2$, 
exactement une
  ou trois des copies
  $s_{\varepsilon_1,\varepsilon_2}(e_1)$,
  $s_{\varepsilon_1,\varepsilon_2}(e_2)$ et
  $s_{\varepsilon_1,\varepsilon_2}(e_3)$ sont effacées.
\end{itemize}

\vspace{1ex}
Appelons le résultat une
\textit{courbe tropicale réelle}. Par exemple si $C$ est une droite
tropicale, il est possible d'effacer six des copies
symétriques des arêtes de $C$ suivant ces deux règles afin d'obtenir
la droite tropicale réelle représentée à la figure \ref{patch
  dte}c.
Si cette courbe tropicale réelle
  n'est certes pas une droite, elle réalise quand même
le même arrangement 
qu'une droite classique dans $\RR^2$ (voir figure \ref{patch dte}d)!

 Il ne s'agit pas là d'un hasard, mais d'un théorème. 

\begin{thm}[O. Viro]\label{viro}
Toute courbe tropicale réelle de degré $d$ réalise le même arrangement
qu'une  courbe algébrique réelle de degré $d$.
\end{thm}

Insistons sur la beauté et la profondeur d'un tel énoncé! En
effet, une courbe tropicale réelle se construit suivant des règles
du jeu combinatoires, et cela ressemble fort à un tour de magie
d'affirmer qu'elle puisse avoir un rapport 
quelconque  avec une courbe algébrique
réelle! Nous ne le ferons pas ici, mais le théorème de Viro nous permet même  
de déterminer l'équation d'une courbe algébrique réelle réalisant
le même arrangement qu'une courbe tropicale réelle donnée.
Utilisons maintenant ce théorème  pour montrer l'existence de deux
courbes algébriques 
réelles, une  de degré 3 et l'autre de degré 6.

\begin{figure}[h]
\begin{center}
\begin{tabular}{ccccccc}
\includegraphics[width=3cm, angle=0]{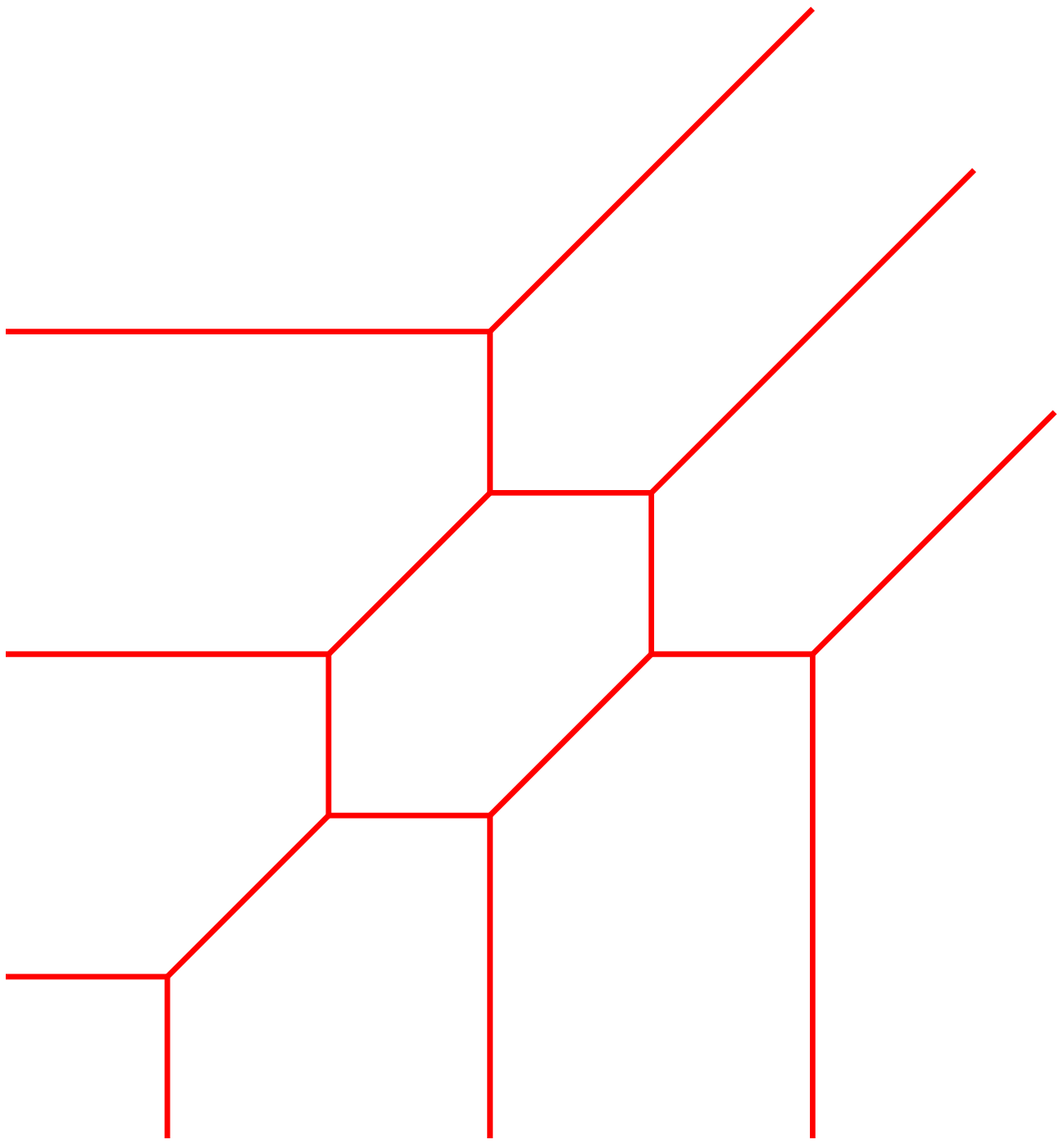}&\hspace{3ex} &
\includegraphics[width=3cm, angle=0]{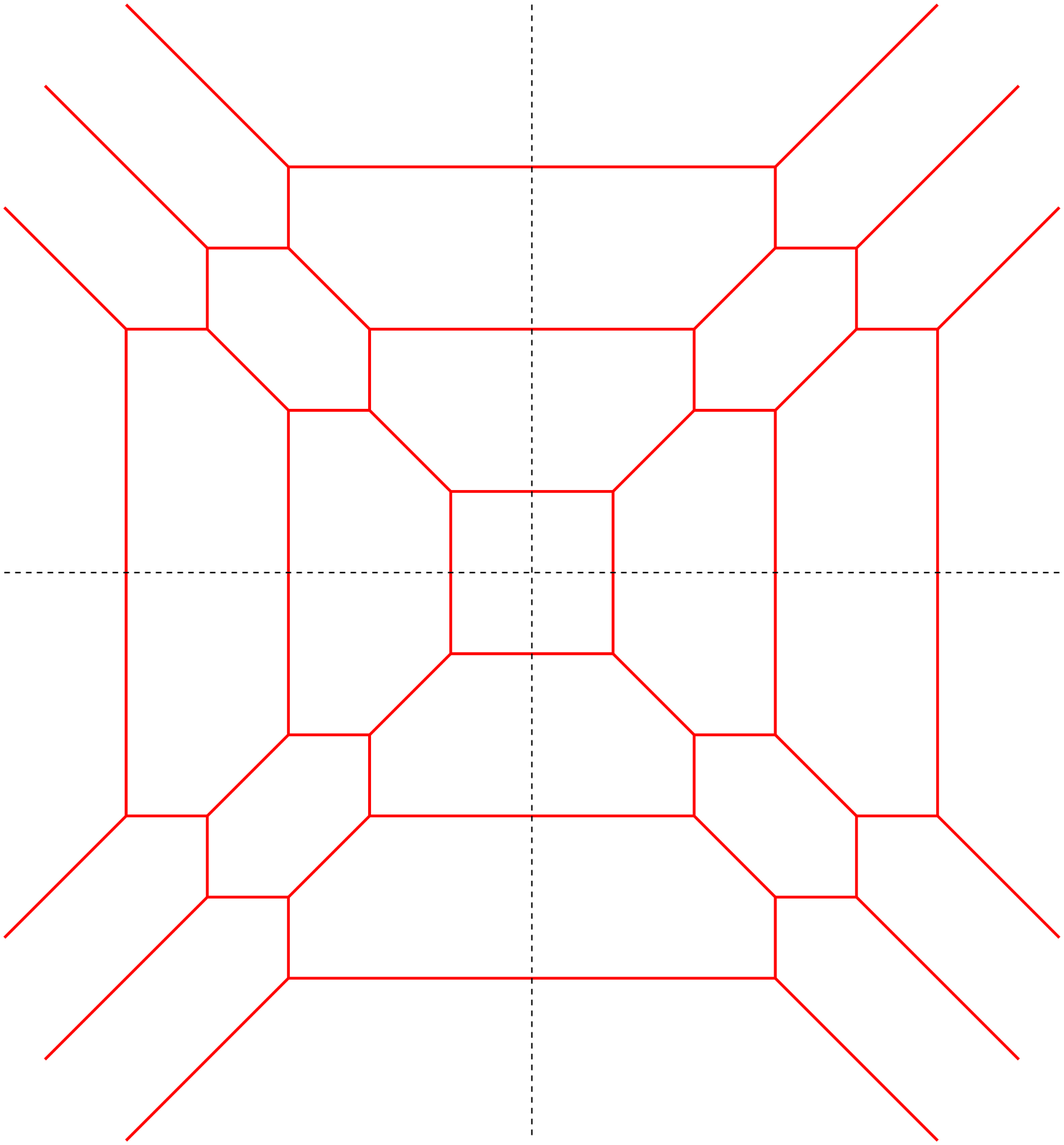}&\hspace{3ex} &
\includegraphics[width=3cm, angle=0]{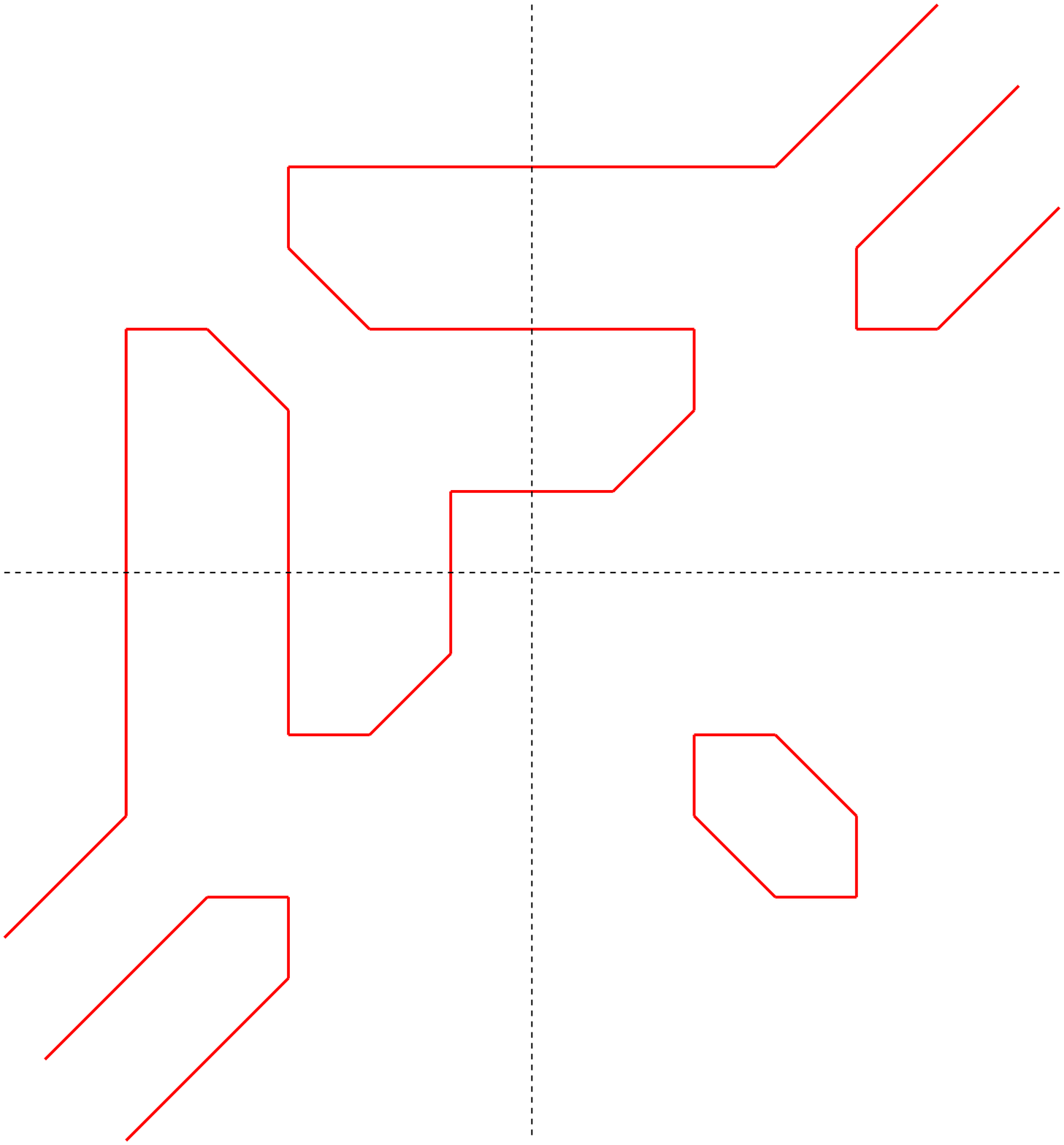}&\hspace{3ex} &
\includegraphics[width=3cm, angle=0]{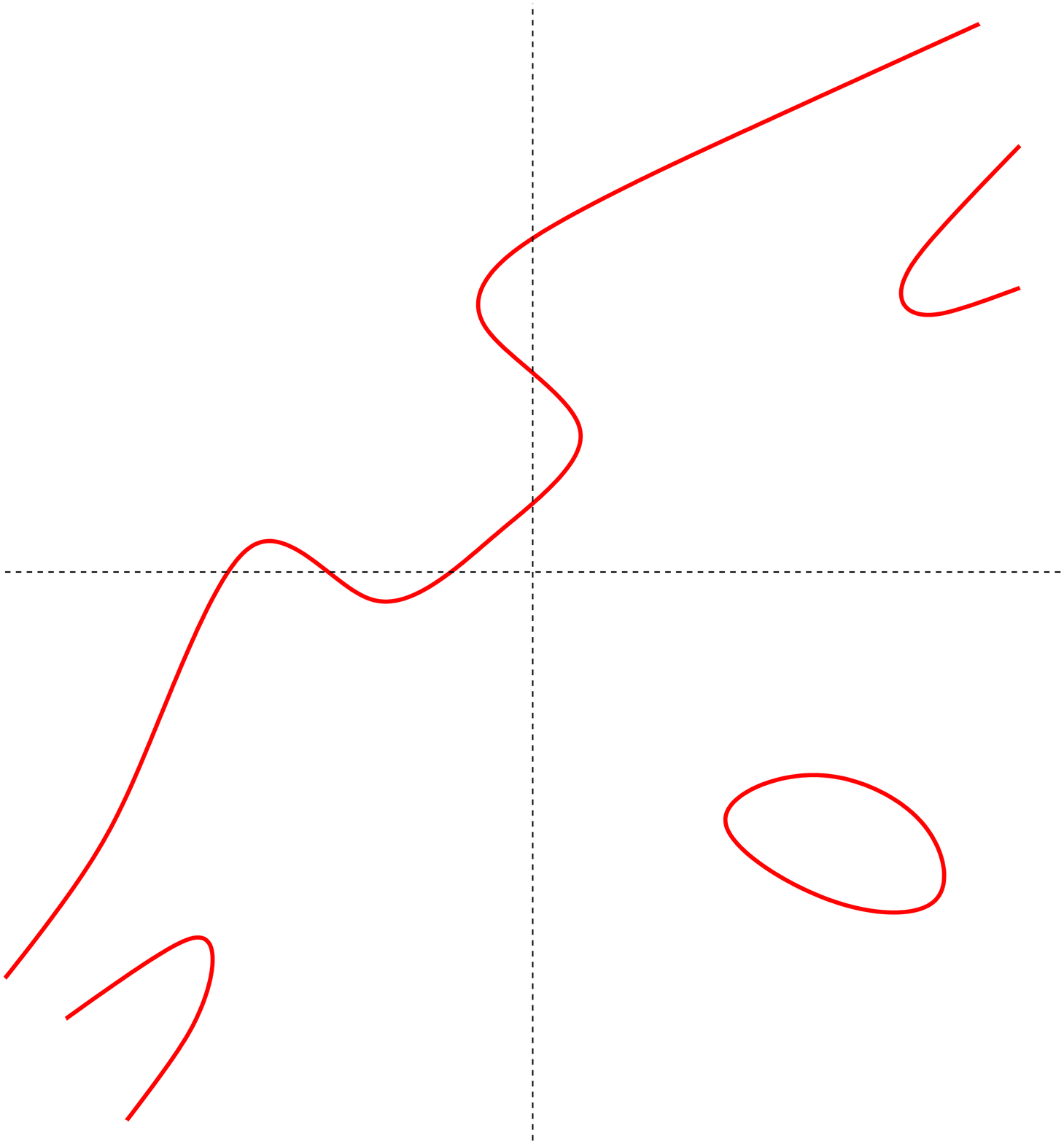}
\\ \\ a) && b) && c) &&d)
\end{tabular}
\end{center}
\caption{Patchwork d'une cubique}
\label{Cub}
\end{figure}

Tout d'abord, considérons la courbe 
tropicale de degré 3 représentée à la figure \ref{Cub}a. 
Pour un choix convenable d'arêtes à effacer, les  figures \ref{Cub}b
et c représentent les deux étapes de la 
 procédure
de patchwork. On a donc prouvé
l'existence d'une courbe algébrique réelle de degré 3 ressemblant 
au dessin de la figure  \ref{Cub}d.

\begin{figure}[h]
\begin{center}
\begin{tabular}{ccc}
\includegraphics[width=5cm, angle=0]{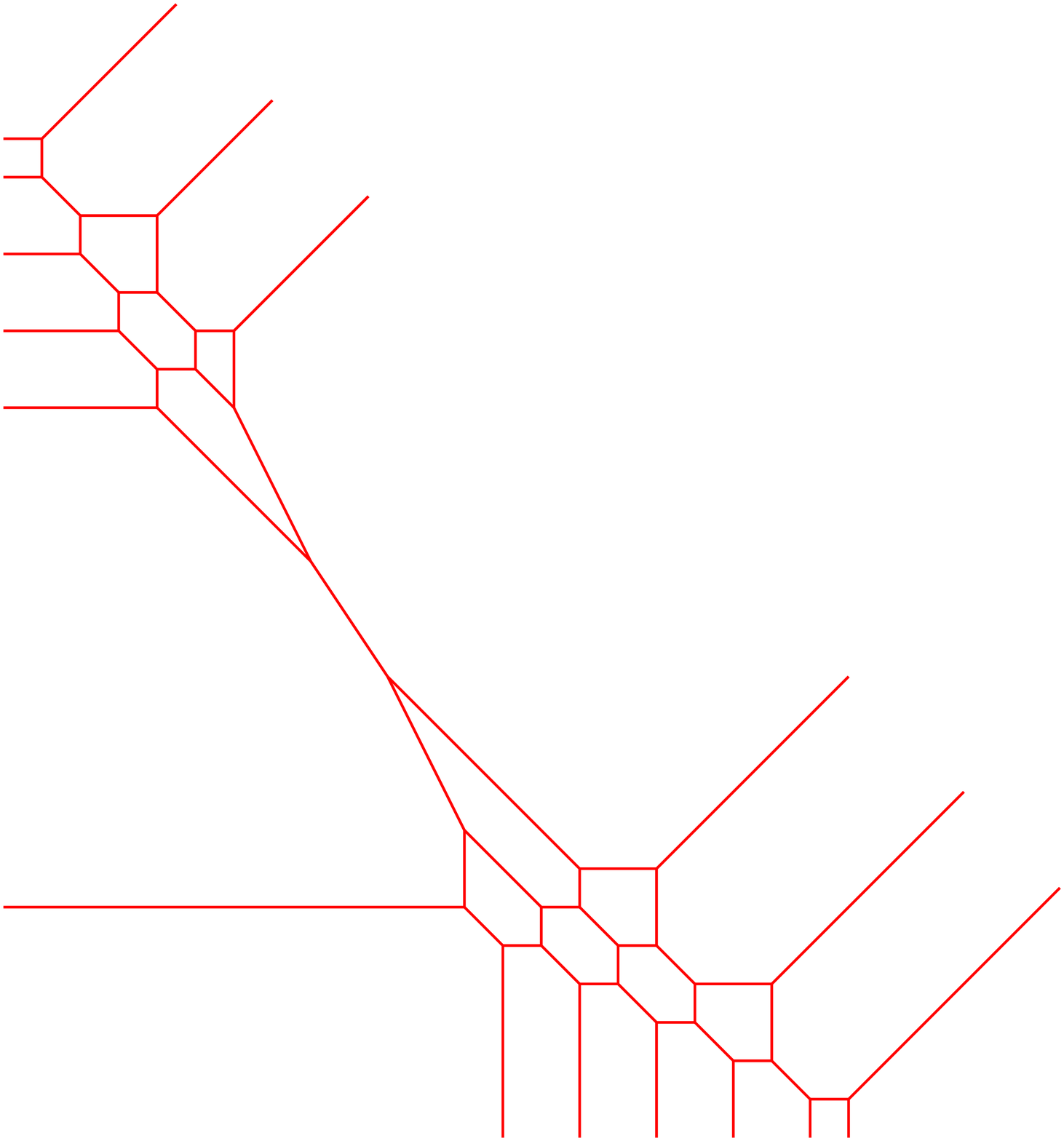}&
\includegraphics[width=5cm, angle=0]{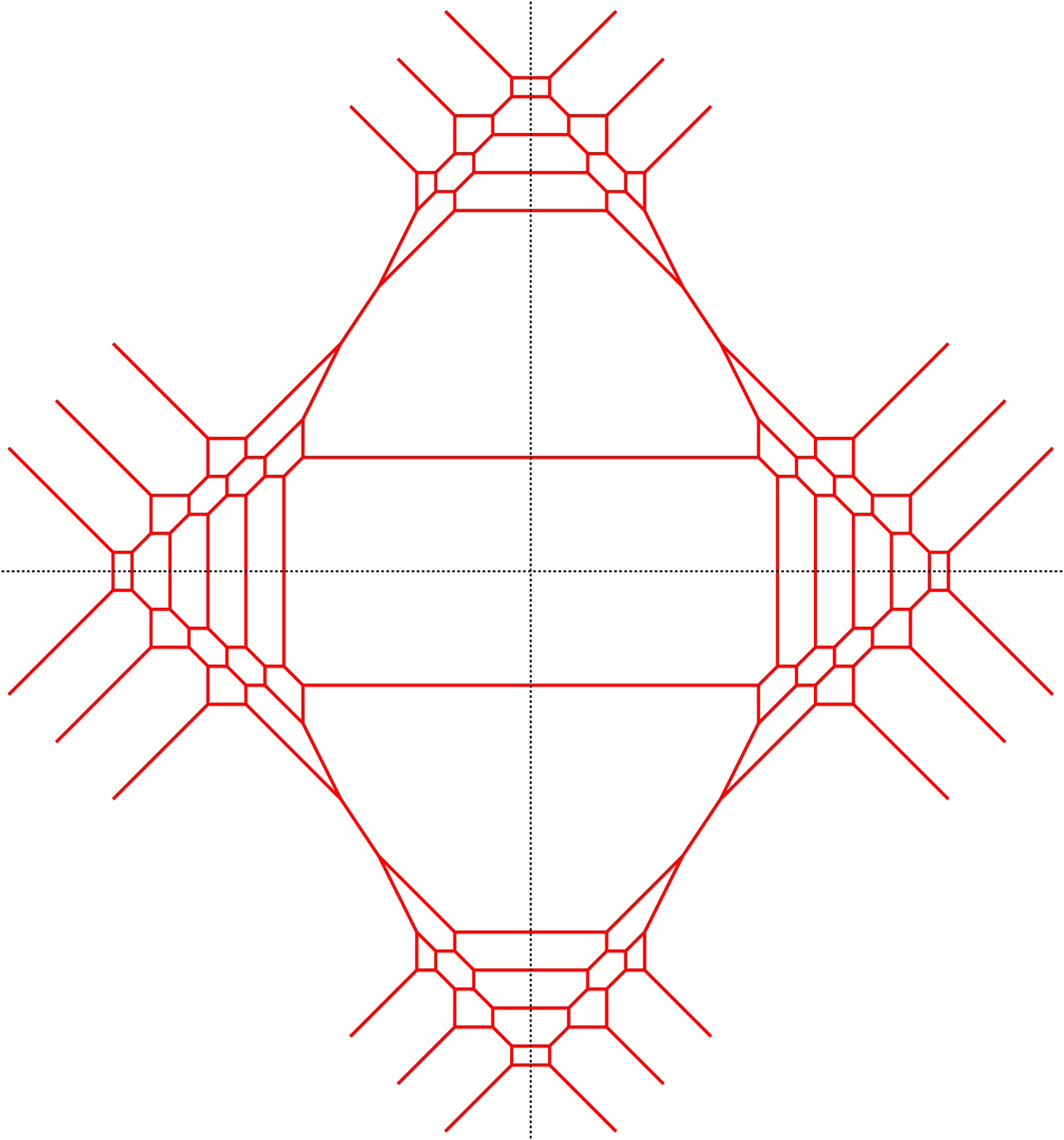}&
\includegraphics[width=5cm, angle=0]{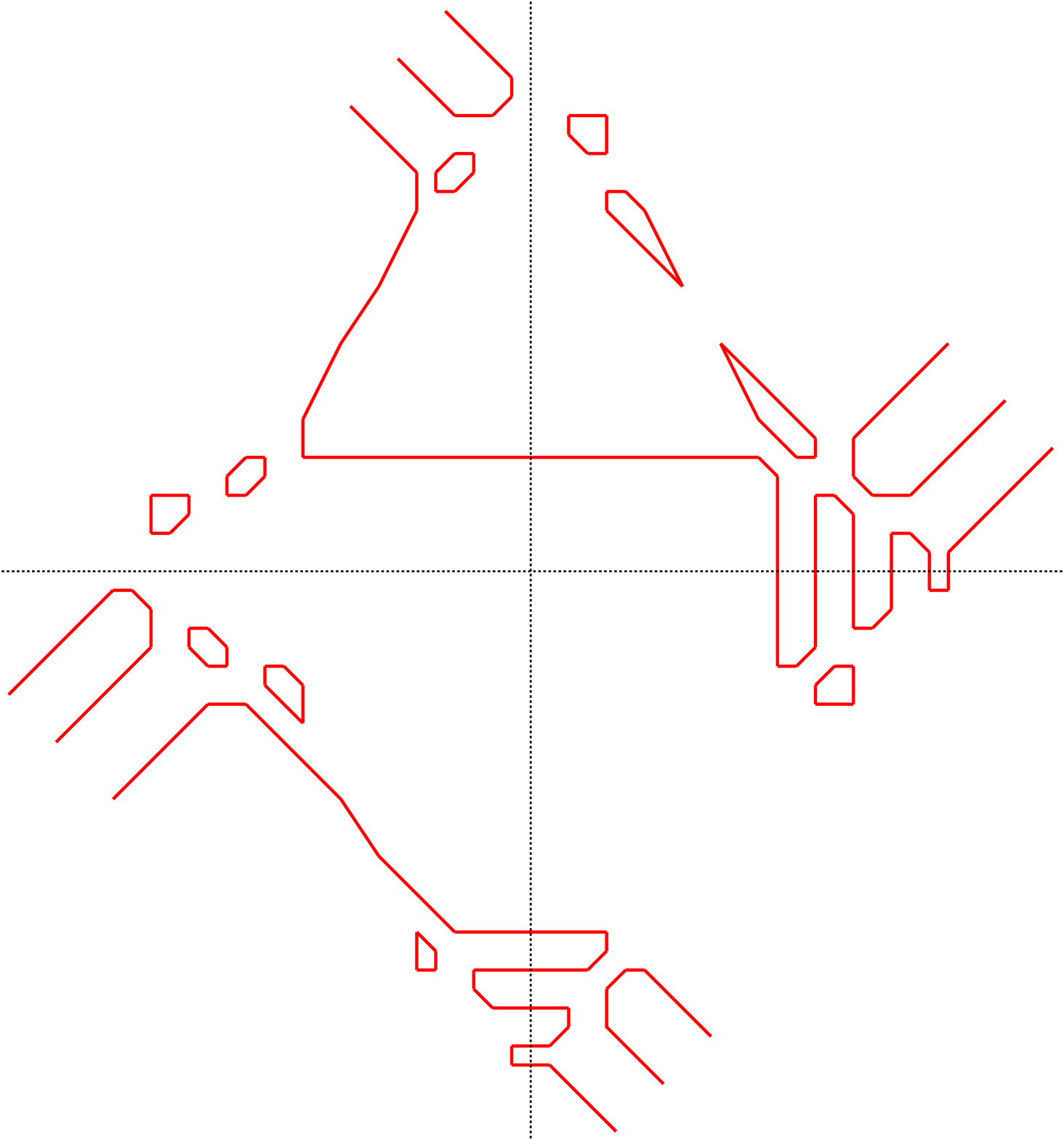}
\\ \\ a) & b) & c)

\end{tabular}
\end{center}
\caption{Courbe de Gudkov}
\label{Gud}
\end{figure}

Pour finir, considérons la courbe 
tropicale de degré 6 représentée à la figure \ref{Gud}a.
Pour un choix convenable d'arêtes à effacer, la procédure
de patchwork donne la courbe de la figure \ref{Gud}c. Une courbe
algébrique réelle de degré 6 réalisant le même arrangement que cette courbe
tropicale réelle a 
été initialement construite par Gudkov, par des moyens beaucoup plus
compliqués,  dans les années 60.
Pour la petite histoire, Hilbert affirmait en 1900 qu'une telle courbe
ne pouvait pas exister...

\subsection{Amibes}
Si la déquantification d'une droite est l'idée qui sous-tend le
patchwork dans toute sa généralité, la preuve du théorème de Viro
est un peu plus technique à écrire
rigoureusement. Contentons nous d'en esquisser les contours.

Tout d'abord, le corps $\RR$ n'étant pas algébriquement clos, il nous
faut  travailler non pas avec des courbes algébriques réelles,
mais plus généralement avec des \textit{courbes algébriques
  complexes}, c'est-à-dire les sous-ensembles de
$(\CC^*)^2$ définis  par
une équation de la forme $P(x,y)=0$ où $P(x,y)$ est un polynôme à coefficients
complexes (qui peuvent donc être réels). 
Pour $t$ un nombre réel positif, on définit l'application
$\text{Log}_t$ sur $(\CC^*)^2$ par: 
$$\begin{array}{cccc}
\text{Log}_t & (\CC^*)^2& \longrightarrow & \RR^2
\\ & (z,w) &\longmapsto  & (\log_t |z|,\log_t |w|)
\end{array} $$

L'image d'une courbe algébrique d'équation  $P(x,y)=0$  par l'application
$\text{Log}_t$, notée $\A_t(P)$, est appelée l'\textit{amibe en base
  $t$} de la courbe. 
Le théorème suivant fournit
un lien fondamental entre la géométrie
algébrique classique et la géométrie tropicale:  toute courbe
tropicale est limite d'amibes de courbes algébriques complexes. 

\begin{thm}[G. Mikhalkin, H. Rullg\aa rd]
Soit $P_{\infty}(x,y)=\tg \sum_{i,j}a_{i,j}x^iy^j \td$ un polynôme
tropical, et soit $\alpha_{i,j}$ un nombre complexe non nul pour chaque
coefficient $a_{i,j}$ différent de $-\infty$. Pour tout $t>0$, on
définit le polynôme 
complexe $P_t(x,y)=\sum_{i,j}\alpha_{i,j}t^{-a_{i,j}}x^iy^j$. Alors
l'amibe $\A_t(P_t)$ converge
vers la courbe tropicale définie par
$P_\infty(x,y)$ lorsque $t$ tend vers l'infini.
\end{thm}
La déquantification de la droite vue à la partie \ref{deqdte} est un
cas particulier de cet énoncé: l'amibe en base $t$ de la droite
d'équation $t^0x-t^0y+t^01=0$ converge vers la droite tropicale
définie par $\tg
0x+0y+0\td$. 
On déduit le théorème de Viro  du théorème précédent en remarquant,
entre autre, que si les $\alpha_{i,j}$ sont des nombres réels, alors
les courbes définies par les polynômes $P_t(x,y)$ sont des courbes
algébriques réelles.

\subsection{Exercices}
\begin{exo}
\begin{enumerate}
\item Construire une courbe tropicale réelle de degré $2$
  réalisant le même arrangement qu'une
  hyperbole dans $\RR^2$. Même question avec une parabole. 
Peut-on construire une courbe tropicale réelle réalisant le même arrangement
qu'une ellipse?
\item À l'aide du patchwork, montrer qu'il existe une courbe
  algébrique réelle de degré 4 réalisant l'arrangement de la figure
  \ref{quartic}b. On pourra s'inspirer de la construction illustrée à la figure
  \ref{Cub}.
\item Montrer que pour tout degré $d$, il existe une courbe algébrique
  réelle plane avec $\frac{d(d-1)+2}{2}$ composantes connexes.
\end{enumerate}
\end{exo}

\section{Références}\label{ref}

Afin de ne pas noyer le lecteur dans un flot de références plus ou
moins accessibles, nous ne renvoyons que vers des textes
d'introduction à la géométrie tropicale et ses applications. 
Pour avoir des références plus spécialisées, on pourra se reporter aux
références des textes cités. Attention, certains auteurs préfèrent
utiliser le minimum au lieu du maximum dans l'algèbre tropicale!

Les introductions à la géométrie tropicale  \cite{BIT} et \cite{St5}
s'adressent à des lecteurs ayant 
un bagage mathématique minimum. Les
lecteurs  plus expérimentés pourront également lire les ouvrages
\cite{St2}, \cite{Mik9} ou \cite{Gath1}. Pour les géomètres confirmés,
nous conseillons les états de l'art \cite{Mik8} et \cite{Mik3}.

Pour en savoir plus sur le 16ème problème de Hilbert, le patchwork, la
déquantification de Maslov et 
les amibes de courbes algébriques, nous renvoyons aux textes
\cite{V9}, \cite{V10}, \cite{IV2}
et \cite{Mik8} ainsi qu'au site Web \cite{Voueb}.

Pour terminer cette introduction à la géométrie tropicale, précisons que
cette dernière s'applique avec succès dans de nombreux domaines autres
que le 16ème problème de Hilbert. Citons par exemple la géométrie
énumérative, la combinatoire, la symétrie miroir, la biologie
mathématique...

\bibliographystyle{alpha}
\bibliography{../../FrBiblio.bib}

\end{document}